\UseRawInputEncoding

\documentclass[journal]{IEEEtran}
\usepackage{amssymb}
\usepackage{amsmath}
\usepackage{pifont,calc}
\usepackage{cite}

\usepackage{graphicx}
\ifCLASSINFOpdf
   \usepackage[pdftex]{graphicx}
\else
 \fi
\usepackage{array}

\usepackage{subeqnarray}
\usepackage{cases} % to combine equations together
\usepackage{bm} % to make xila zimu hei ti
\usepackage{subfigure}
\usepackage{mathrsfs}
\usepackage{url}
\usepackage[ruled]{algorithm2e}
\ifCLASSINFOpdf
\else
\fi
\begin{document}
\title{Efficient and Stable Algorithms to Extend Greville's Method to Partitioned Matrices Based on
Inverse Cholesky Factorization}

\author{Hufei~Zhu
\thanks{H. Zhu is with College of Computer Science and Software, Shenzhen University, Shenzhen 518060, China (e-mail: zhuhufei@aliyun.com).}}


\markboth{Journal of \LaTeX\ Class Files,~Vol.~14, No.~8, August~2015}%
{Shell \MakeLowercase{\textit{et al.}}: Bare Demo of IEEEtran.cls for IEEE Journals}

%\begin{equation}\label{equ1}
%
% \end{equation}
%
% \begin{equation}\label{equ1}
%
% \end{equation}
%
% \begin{equation}\label{equ1}
%
% \end{equation}



\maketitle

\begin{abstract}
Greville's method has been utilized in (Broad Learning System) BLS to propose an effective and efficient
incremental learning system without
retraining the whole network from the beginning.
For a column-partitioned matrix
where the second part consists of  $p$ columns,
Greville's method requires  $p$ iterations to compute the pseudoinverse of the whole matrix from the pseudoinverse of the first part.
The incremental algorithms in BLS extend Greville's method
   to compute  the pseudoinverse of the whole matrix from the pseudoinverse of the first part
    by just $1$ iteration,
   which have neglected some possible cases, and need further improvements in efficiency and numerical stability.
   In this paper, we propose an efficient and numerical stable algorithm from Greville's method,
 to compute  the pseudoinverse of the whole matrix from the pseudoinverse of the first part by just $1$ iteration,
where all possible cases are considered, and  the recently proposed inverse Cholesky factorization
 can be applied to further reduce the computational complexity. Finally, we give the whole algorithm for column-partitioned matrices in BLS.
 On the other hand, we also give the proposed algorithm for row-partitioned matrices in BLS.
\end{abstract}


% \fi
\IEEEpeerreviewmaketitle

\section{Introduction}

In \cite{BL_trans_paper}, the pseudoinverse of a partitioned
matrix, i.e., Greville's method \cite{cite_general_inv_book}, has been utilized to propose Broad Learning System (BLS), an effective and efficient
incremental learning system without
retraining the whole network from the beginning.
For a column-partitioned matrix
%\begin{equation}\label{Abig_original}
$\mathbf{A}_{{}}^{m+1}=\left[ {{\mathbf{A}}^{m}}|\mathbf{H}_{m+1}^{{}} \right]$
%\end{equation}
where
%the $m \times p$
 $\mathbf{H}_{m+1}^{{}}$ has $p$ columns,
Greville's method \cite{cite_general_inv_book}  needs $p$ iterations to
compute ${{\left( \mathbf{A}_{{}}^{m+1} \right)}^{+ }}$ from ${{\left( \mathbf{A}_{{}}^{m} \right)}^{+ }}$,
where ${{\mathbf{A}}^{+ }}$ denotes the pseudoinverse of the matrix $\mathbf{A}$.  However, the incremental algorithms in \cite{BL_trans_paper}
apply Greville's method
   to compute ${{\left( \mathbf{A}_{{}}^{m+1} \right)}^{+ }}$ from ${{\left( \mathbf{A}_{{}}^{m} \right)}^{+ }}$ by just $1$ iteration,
   which have neglected some possible cases, and need further improvements in efficiency and numerical stability.
Based on Greville's method, we propose an efficient and numerical stable algorithm to compute ${{\left( \mathbf{A}_{{}}^{m+1} \right)}^{+ }}$ from ${{\left( \mathbf{A}_{{}}^{m} \right)}^{+ }}$ by just $1$ iteration,
where all possible cases are considered, and  the efficient inverse Cholesky factorization in \cite{my_inv_chol_paper}
 can be applied to further reduce the computational complexity.
 Moreover, the proposed algorithm is also applied to  row-partitioned matrices in BLS.

\section{The Proposed Improvements for Greville's method Utilized in BLS}

In subsection A, we introduce Greville's method \cite{cite_general_inv_book} and
its application to
column-partitioned matrices in BLS. In subsection B, we deduce three  theorems that will be utilized.
Then in subsection C, we propose the modified Greville's method for BLS, which considers all possible cases,
 and is improved in efficiency and numerical stability.  Finally in subsection D,
 we apply  the efficient inverse Cholesky factorization in \cite{my_inv_chol_paper} to further reduce the computational complexity.

\subsection{Greville's method and its application to column-partitioned matrices in BLS}

As in \cite{BL_trans_paper},
write
 the  $m \times (n+p)$ column-partitioned matrix as
%(1).ÊéÉϵÄԭ'¹«Ê½
\begin{equation}\label{Abig_original}
\mathbf{A}_{{}}^{m+1}=\left[ {{\mathbf{A}}^{m}}|\mathbf{H}_{m+1}^{{}} \right],
\end{equation}
where ${{\mathbf{A}}^{m}}$ is  $m \times n$  and  $\mathbf{H}_{m+1}^{{}}$
is $m \times p$.
Let
\begin{equation}\label{A_k_def12433}
\mathbf{A}_{k}^{m+1}=\left[ {{\mathbf{A}}^{m}}|\mathbf{H}_{m+1}^{:,1:k} \right],
\end{equation}
 where $\mathbf{H}_{m+1}^{:,1:k}$ denotes the first $k$ columns of $\mathbf{H}_{m+1}^{{}}$.  Then
 % we have
\begin{equation}\label{A_k_def_iterative}
\mathbf{A}_{k}^{m+1}=\left[ \mathbf{A}_{k-1}^{m+1}|\mathbf{h}_{:k}^{m+1} \right],
\end{equation}
where $\mathbf{h}_{:k}^{m+1}$ denotes the $k$-th column of $\mathbf{H}_{m+1}^{{}}$. Notice that when $k=0$, $\mathbf{H}_{m+1}^{:,1:k}$ becomes empty and then (\ref{A_k_def12433}) becomes
\begin{equation}\label{A0isAm}
\mathbf{A}_{0}^{m+1}={{\mathbf{A}}^{m}}.
\end{equation}

Greville's method \cite{cite_general_inv_book}   computes ${{\left( \mathbf{A}_{p}^{m+1} \right)}^{+ }}={{\left( \mathbf{A}_{{}}^{m+1} \right)}^{+ }}$ from ${{\left( \mathbf{A}_{0}^{m+1} \right)}^{+ }}={{\left( \mathbf{A}_{{}}^{m} \right)}^{+ }}$ by $p$ iterations.
In the $k$-th iteration ($k=1,2,\cdots, p$), ${{\left( \mathbf{A}_{k}^{m+1} \right)}^{+ }}$  is \cite[Theorem 5.7]{cite_general_inv_book}
\begin{equation}\label{A_1col_inv_book1}{{\left[ \begin{matrix}
   \mathbf{A}_{k-1}^{m+1} & \mathbf{h}_{:k}^{m+1}  \\
\end{matrix} \right]}^{+}}=\left[ \begin{matrix}
   {{(\mathbf{A}_{k-1}^{m+1})}^{+}}-{{{\mathbf{\tilde{d}}}}_{:k}}\mathbf{\tilde{b}}_{:k}^{T}  \\
   \mathbf{\tilde{b}}_{:k}^{T}  \\
\end{matrix} \right],
\end{equation}
where
 \begin{equation}\label{d_col_define4311}
 {{\mathbf{\tilde{d}}}_{:k}}={{(\mathbf{A}_{k-1}^{m+1})}^{+}}\mathbf{h}_{:k}^{m+1},
 \end{equation}
and $\mathbf{\tilde{b}}_{:k}^{T}$ is computed from
 \begin{equation}\label{c_col_define1}
 \mathbf{\tilde{c}}_{:k}^{{}}=\mathbf{h}_{:k}^{m+1}-\mathbf{A}_{k-1}^{m+1}{{\mathbf{\tilde{d}}}_{:k}}
 \end{equation}
by
%\begin{equation}\label{b_col_def1ab}
%\mathbf{\tilde{b}}_{:k}^{T}=\left\{ \begin{matrix}
%   \mathbf{\tilde{c}}_{:k}^{+} & i{{f}^{{}}}\mathbf{\tilde{c}}_{:k}^{{}}\ne \mathbf{0},  \\
%   {{(1+\mathbf{\tilde{d}}_{:k}^{T}{{{\mathbf{\tilde{d}}}}_{:k}})}^{-1}}\mathbf{\tilde{d}}_{:k}^{T}{{(\mathbf{A}_{k-1}^{m+1})}^{+}} & i{{f}^{{}}}\mathbf{\tilde{c}}_{:k}^{{}}=\mathbf{0}.  \\
%\end{matrix} \right.
%\end{equation}
%\begin{subnumcases}{\label{b_col_def1ab}}
%{{\bf{\tilde b}}_{:k}^T = {\bf{\tilde c}}_{:k}^ + } \quad if \ {\bf{\tilde c}}_{:k}^{} \ne {\bf{0}} &  \label{b_col_def1a}\\
%{{\bf{\tilde b}}_{:k}^T = {{(1 + {\bf{\tilde d}}_{:k}^T{{{\bf{\tilde d}}}_{:k}})}^{ - 1}}{\bf{\tilde d}}_{:k}^T{{({\bf{A}}_{k - 1}^{m + 1})}^ + }} \quad if \ {\bf{\tilde c}}_{:k}^{} = {\bf{0}}. & \label{b_col_def1b}
%\end{subnumcases}
\begin{subequations}{\label{b_col_def1ab}}
 \begin{numcases}
{ {\bf{\tilde b}}_{:k}^T =   }
{ {\bf{\tilde c}}_{:k}^ + } \quad \quad \quad \quad \quad \quad \quad \quad \quad \quad \quad \   if \ {\bf{\tilde c}}_{:k}^{} \ne {\bf{0}}   &  \label{b_col_def1a} \\
{{{(1 + {\bf{\tilde d}}_{:k}^T{{{\bf{\tilde d}}}_{:k}})}^{ - 1}}{\bf{\tilde d}}_{:k}^T{{({\bf{A}}_{k - 1}^{m + 1})}^ + }} \quad if \ {\bf{\tilde c}}_{:k}^{} = {\bf{0}}.  &  \label{b_col_def1b}
\end{numcases}
\end{subequations}
Here
${{\mathbf{A}}^{+ }}$ is   the unique Moore-Penrose generalized inverse (i.e., the pseudoinverse) that
% of the matrix $\mathbf{A}$, which
satisfies \cite{cite_general_inv_book}
\begin{subnumcases}{\label{MoorePenroseDef121all}}
\mathbf{A}{{\mathbf{A}}^{+ }}\mathbf{A}=\mathbf{A} &  \label{MoorePenroseDef121a}\\
{{\mathbf{A}}^{+ }}\mathbf{A}{{\mathbf{A}}^{+ }}={{\mathbf{A}}^{+ }}  &  \label{MoorePenroseDef121b}\\
{{(\mathbf{A}{{\mathbf{A}}^{+ }})}^{T}}=\mathbf{A}{{\mathbf{A}}^{+ }}   &  \label{MoorePenroseDef121c}\\
{{({{\mathbf{A}}^{+ }}\mathbf{A})}^{T}}={{\mathbf{A}}^{+ }}\mathbf{A}. & \label{MoorePenroseDef121d}
\end{subnumcases}
%Greville's method \cite{cite_general_inv_book}  computes ${{\left( \mathbf{A}_{p}^{m+1} \right)}^{+ }}={{\left( \mathbf{A}_{{}}^{m+1} \right)}^{+ }}$ from ${{\left( \mathbf{A}_{0}^{m+1} \right)}^{+ }}={{\left( \mathbf{A}_{{}}^{m} \right)}^{+ }}$ by $p$ iterations.

%(2).´ÓÁÐÏòÁ¿Íƹ㵽¾ØÕó
In \cite{BL_trans_paper}, the column vector $\mathbf{h}_{:k}^{m+1}$  in
(\ref{A_1col_inv_book1})
 is extended to the matrix ${{\mathbf{H}}_{m+1}}$  with $p$ columns, and correspondingly
 (\ref{A_1col_inv_book1}),
(\ref{d_col_define4311}),
(\ref{c_col_define1})
and
(\ref{b_col_def1ab})
 become
 \begin{equation}\label{A_1col_inv_book1matrix}
 {{({{\mathbf{A}}^{m+1}})}^{+}}={{\left[ {{\mathbf{A}}^{m}}|{{\mathbf{H}}_{m+1}} \right]}^{+}}=\left[ \begin{matrix}
   {{({{\mathbf{A}}^{m}})}^{+}}-\mathbf{D}{{\mathbf{B}}^{T}}  \\
   {{\mathbf{B}}^{T}}  \\
\end{matrix} \right],
\end{equation}
\begin{equation}\label{D_matrix_def_111}
\mathbf{D} = {{({{\mathbf{A}}^{m}})}^{+}}{{\mathbf{H}}_{m+1}}\text{ },
\end{equation}
\begin{equation}\label{C_matrix_def_111}
\mathbf{C}={{\mathbf{H}}_{m+1}}-{{\mathbf{A}}^{m}}\mathbf{D},
\end{equation}
and
%\begin{equation}\label{}{{\mathbf{B}}^{T}}=\left\{ \begin{matrix}
%   {{\mathbf{C}}^{+}} & i{{f}^{{}}}\mathbf{C}\ne 0,  \\
%   {{(1+{{\mathbf{D}}^{T}}\mathbf{D})}^{-1}}{{\mathbf{B}}^{T}}{{({{\mathbf{A}}^{m}})}^{+}} & i{{f}^{{}}}\mathbf{C}=0,  \\
%\end{matrix} \right.\end{equation}         (B_Matrix_def1a,b)
%\begin{subnumcases}{\label{B_Matrix_def1ab}}
%{{{\bf{B}}^T} = {{\bf{C}}^ + } \quad if \ {\bf{C}} \ne 0} &  \label{B_Matrix_def1a}\\
%{{{\bf{B}}^T} = {{(1 + {{\bf{D}}^T}{\bf{D}})}^{ - 1}}{{\bf{D}}^T}{{({{\bf{A}}^m})}^ + } \quad if \ {\bf{C}} = 0}, & \label{B_Matrix_def1b}
%\end{subnumcases}
\begin{subequations}{\label{B_Matrix_def1ab}}
 \begin{numcases}
{  {{\bf{B}}^T} = }
{ {{\bf{C}}^ + } \quad \quad \quad \quad \quad \quad \quad \quad \quad \quad   if \ {\bf{C}} \ne 0}   &  \label{B_Matrix_def1a} \\
{{{({\bf{I}} + {{\bf{D}}^T}{\bf{D}} )}^{ - 1}}{{\bf{D}}^T}{{({{\bf{A}}^m})}^ + } \quad if \ {\bf{C}} = 0},  &  \label{B_Matrix_def1b}
\end{numcases}
\end{subequations}
respectively.
%(3)ÈÎÒâc_waveΪ0¾ÍÐèÒªÓÃd¡£
%Whether ${\bf{\tilde c}}_{:k}^{} \ne {\bf{0}}$ and ${\bf{C}} \ne 0$  are checked in
%(\ref{b_col_def1ab}) and (\ref{B_Matrix_def1ab}), respectively,
%but
%%Now let us compare
% $\mathbf{\tilde{c}}_{:k}^{{}}$ is usually different from  $\mathbf{c}_{:k}^{{}}$, the $k$-th column of ${\bf{C}}$, as shown in the following.

${\bf{C}} \ne 0$ is required in (\ref{B_Matrix_def1a}), while  ${\bf{\tilde c}}_{:k}^{} \ne {\bf{0}}$ is required in (\ref{b_col_def1a}).
Denote the $k$-th column of ${\bf{C}}$ as $\mathbf{c}_{:k}^{{}}$. In the next paragraph we will show that usually $\mathbf{c}_{:k}^{{}}$ is different from
 $\mathbf{\tilde{c}}_{:k}^{{}}$.
  Then it can be easily seen that (\ref{b_col_def1a}) cannot be  extended to (\ref{B_Matrix_def1a}),
 since $\mathbf{c}_{:k}^{{}} \ne {\bf{0}}$ for each $k$  ($1 \le k \le p$) in (\ref{B_Matrix_def1ab})
 %does not mean
 can not ensure
  ${\bf{\tilde c}}_{:k}^{} \ne {\bf{0}}$ in (\ref{b_col_def1a}),
 and actually ${\bf{C}} \ne 0$ only means at least one $\mathbf{c}_{:k}^{{}} \ne {\bf{0}}$.

 %as shown in the next paragraph.
%in what follows,
To show that usually $\mathbf{c}_{:k}^{{}}$ is different from
 $\mathbf{\tilde{c}}_{:k}^{{}}$,
 %$\mathbf{c}_{:k}^{{}} \ne \mathbf{\tilde{c}}_{:k}^{{}}$,
substitute  (\ref{d_col_define4311})
into  (\ref{c_col_define1})  to obtain
\begin{equation}\label{c_wave_k_final_compute123}
\mathbf{\tilde{c}}_{:k}^{{}}=\mathbf{h}_{:k}^{m+1}-\mathbf{A}_{k-1}^{m+1}{{(\mathbf{A}_{k-1}^{m+1})}^{+}}\mathbf{h}_{:k}^{m+1},
\end{equation}
%On the other hand,
and substitute (\ref{D_matrix_def_111}) into (\ref{C_matrix_def_111}) to obtain
\begin{equation}\label{C_matrix_def_onlyH121}
\mathbf{C}={{\mathbf{H}}_{m+1}}-{{\mathbf{A}}^{m}}{{({{\mathbf{A}}^{m}})}^{+}}{{\mathbf{H}}_{m+1}}\text{ },
\end{equation}
from which we can deduce %that the $k$-th column of ${\bf{C}}$ is
\begin{equation}\label{c_k_final_compute123}
\mathbf{c}_{:k}^{{}}=\mathbf{h}_{:k}^{m+1}-\mathbf{A}_{{}}^{m}{{(\mathbf{A}_{{}}^{m})}^{+}}\mathbf{h}_{:k}^{m+1}.
\end{equation}
From  (\ref{c_k_final_compute123}) and (\ref{c_wave_k_final_compute123}), it can be seen that only
\begin{equation}\label{c1_wave_c1_equal}
\mathbf{\tilde{c}}_{:1}^{{}}\text{=}\mathbf{c}_{:1}^{{}},
\end{equation}
and $for \quad k=2,3,\cdots, p$, usually
\begin{equation}\label{c2k_wave_c2k_NotEqual}
\mathbf{\tilde{c}}_{:k}^{{}}\ne \mathbf{c}_{:k}^{{}}
\end{equation}
since usually
%$\mathbf{A}_{{}}^{m}{{(\mathbf{A}_{{}}^{m})}^{+}}\ne \mathbf{A}_{k-1}^{m+1}{{(\mathbf{A}_{k-1}^{m+1})}^{+}}$.
\begin{equation}\label{AinvAnotEQUAkinvAk}
\mathbf{A}_{{}}^{m}{{(\mathbf{A}_{{}}^{m})}^{+}}\ne \mathbf{A}_{k-1}^{m+1}{{(\mathbf{A}_{k-1}^{m+1})}^{+}}.
 \end{equation}

 \subsection{Three Theorems about $\mathbf{\tilde{c}}_{:k}$}
 %Relevant Theorems

 In this subsection, we deduce three  theorems relevant to $\mathbf{\tilde{c}}_{:k}$.

%\subsubsection{Theorem 1}

%(3.5)ÏÂÃæverifyÒÔÇ°c_wave¶¼ÊÇÁãʱ£¬c_wave = CµÄÁÐ
Firstly about (\ref{AinvAnotEQUAkinvAk}),
 we have
% Let $\Re (\mathbf{A})$ denote the range of any $\mathbf{A}\in {{\mathbb{R}}^{N\times k}}$, i.e.,
%\begin{equation}\label{general_inverse_book}
%\Re (\mathbf{A})\text{=}\left\{ \mathbf{y}\in {{\mathbb{R}}^{N}}:\mathbf{y}=\mathbf{Ax}_{{}}^{{}}for_{{}}^{{}}some_{{}}^{{}}\mathbf{x}\in {{\mathbb{R}}^{k}} \right\}.
%\end{equation}
%Since ${{\mathbf{\tilde{c}}}_{:k}}=0$ is equivalent to  $\Re (\mathbf{A}_{k}^{m+1})=\Re (\mathbf{A}_{k-1}^{m+1})$ \cite{cite_general_inv_book},

 \underline{\textbf{Theorem 1}}. If $\mathbf{\tilde{c}}_{:k}^{{}}=\mathbf{0}$, then
\begin{equation}\label{A_equal_c0123}
\mathbf{A}_{k}^{m+1}{{(\mathbf{A}_{k}^{m+1})}^{+}}=\mathbf{A}_{k-1}^{m+1}{{(\mathbf{A}_{k-1}^{m+1})}^{+}}.
\end{equation}

Proof. Applying (\ref{A_1col_inv_book1}) to obtain
 $\mathbf{A}_{k}^{m+1}{{(\mathbf{A}_{k}^{m+1})}^{+}}=\left[ \begin{matrix}
   \mathbf{A}_{k-1}^{m+1} & \mathbf{h}_{:k}^{m+1}  \\
\end{matrix} \right]\left[ \begin{matrix}
   {{(\mathbf{A}_{k-1}^{m+1})}^{+}}-{{{\mathbf{\tilde{d}}}}_{:k}}\mathbf{\tilde{b}}_{:k}^{T}  \\
   \mathbf{\tilde{b}}_{:k}^{T}  \\
\end{matrix} \right]$, i.e.,
\begin{multline}\label{A_equal_c02prove}
\mathbf{A}_{k}^{m+1}{{(\mathbf{A}_{k}^{m+1})}^{+}}=\mathbf{A}_{k-1}^{m+1}{{(\mathbf{A}_{k-1}^{m+1})}^{+}} \\
-\mathbf{A}_{k-1}^{m+1}{{\mathbf{\tilde{d}}}_{:k}}
\mathbf{\tilde{b}}_{:k}^{T}+\mathbf{h}_{:k}^{m+1}\mathbf{\tilde{b}}_{:k}^{T}.
\end{multline}
Now we only need to verify that the last two entries in the right side of (\ref{A_equal_c02prove}) satisfy
\begin{equation}\label{A_last2_c02prove}
-\mathbf{A}_{k-1}^{m+1}{{\mathbf{\tilde{d}}}_{:k}}\mathbf{\tilde{b}}_{:k}^{T}+\mathbf{h}_{:k}^{m+1}\mathbf{\tilde{b}}_{:k}^{T}=\mathbf{0}.
\end{equation}
%Let us substitute (\ref{d_col_define4311}) into
%%the left side of
% (\ref{A_last2_c02prove}),  to write it as
%which can be written as
%\begin{equation}\label{2verifyc0Cuse123}
%-\mathbf{A}_{k-1}^{m+1}{{(\mathbf{A}_{k-1}^{m+1})}^{+}}\mathbf{h}_{:k}^{m+1}\mathbf{\tilde{b}}_{:k}^{T}+\mathbf{h}_{:k}^{m+1}\mathbf{\tilde{b}}_{:k}^{T}=\mathbf{0}
%\end{equation}
%by substituting (\ref{d_col_define4311}) into
%%the left side of
% (\ref{A_last2_c02prove}).
Since $\mathbf{\tilde{c}}_{:k}^{{}}=\mathbf{0}$, from (\ref{c_wave_k_final_compute123}) we can deduce
\begin{equation}\label{c0_AAh2h}
\mathbf{A}_{k-1}^{m+1}{{(\mathbf{A}_{k-1}^{m+1})}^{+}}\mathbf{h}_{:k}^{m+1}=\mathbf{h}_{:k}^{m+1},
\end{equation}
into which we substitute  (\ref{d_col_define4311}) to obtain
\begin{equation}\label{c0_AAh2hMay14a92832}
\mathbf{A}_{k-1}^{m+1}  {{\mathbf{\tilde{d}}}_{:k}}=\mathbf{h}_{:k}^{m+1}.
\end{equation}
Then we can substitute (\ref{c0_AAh2hMay14a92832}) into
(\ref{A_last2_c02prove}) to verify (\ref{A_last2_c02prove}).
%which can be substituted into (\ref{2verifyc0Cuse123}) to verify it.
%Finally we can substitute (\ref{c0_AAh2h}) into (\ref{2verifyc0Cuse123}), to verify that (\ref{2verifyc0Cuse123}) is equal to zero, i.e., we have verified (\ref{A_last2_c02prove}).

On the other hand, notice that ${{\mathbf{\tilde{c}}}_{:k}}=0$ is equivalent \cite[the last 3rd and 4th rows in page 166]{cite_general_inv_book} to
 %$\Re (\mathbf{A}_{k}^{m+1})=\Re (\mathbf{A}_{k-1}^{m+1})$ ,
 \begin{equation}\label{RAkm1equ2RAkm1}
\Re (\mathbf{A}_{k}^{m+1})=\Re (\mathbf{A}_{k-1}^{m+1}),
 \end{equation}
where
%$\Re (\mathbf{A})$ is the range of any $\mathbf{A}\in {{\mathbb{R}}^{N\times k}}$
\begin{equation}\label{general_inverse_book}
\Re (\mathbf{A})\text{=}\left\{ \mathbf{y}\in {{\mathbb{R}}^{N}}:\mathbf{y}=\mathbf{Ax} \ for \ some \ \mathbf{x}\in {{\mathbb{R}}^{k}} \right\}
\end{equation}
is the range  \cite{cite_general_inv_book} of any $\mathbf{A}\in {{\mathbb{R}}^{N\times k}}$.
% we can also deduce (\ref{A_equal_c0123}) from (\ref{RAkm1equ2RAkm1}).

%\subsubsection{Theorem 2}

From Theorem 1,
 we derive

\underline{\textbf{Theorem 2}}. If
\begin{equation}\label{c_wave_1_k_all0_1223}
\mathbf{\tilde{c}}_{:k}^{{}}=\mathbf{\tilde{c}}_{:k-1}^{{}}=\cdots =\mathbf{\tilde{c}}_{:1}^{{}}=\mathbf{0},
\end{equation}
then
\begin{equation}\label{c_wave_ c_k _equal}
\mathbf{\tilde{c}}_{:k+1}^{{}}=\mathbf{c}_{:k+1}^{{}}.
\end{equation}

Proof.  Apply (\ref{A_equal_c0123})  iteratively to obtain
%\begin{equation}\label{AAk0Equal123}
%\mathbf{A}_{k}^{m+1}{{(\mathbf{A}_{k}^{m+1})}^{+}}=\mathbf{A}_{k-1}^{m+1}{{(\mathbf{A}_{k-1}^{m+1})}^{+}}=\mathbf{A}_{k-2}^{m+1}{{(\mathbf{A}_{k-2}^{m+1})}^{+}}=\cdots
%=\mathbf{A}_{0}^{m+1}{{(\mathbf{A}_{0}^{m+1})}^{+}},
%\end{equation}
$\mathbf{A}_{k}^{m+1}{{(\mathbf{A}_{k}^{m+1})}^{+}}=\mathbf{A}_{k-1}^{m+1}{{(\mathbf{A}_{k-1}^{m+1})}^{+}}=\mathbf{A}_{k-2}^{m+1}{{(\mathbf{A}_{k-2}^{m+1})}^{+}}=\cdots
=\mathbf{A}_{0}^{m+1}{{(\mathbf{A}_{0}^{m+1})}^{+}}$,
which can be substituted into (\ref{c_wave_k_final_compute123}) to deduce
\begin{equation}\label{c_wave_k1_A0_123}
\mathbf{\tilde{c}}_{:k+1}^{{}}=\mathbf{h}_{:k+1}^{m+1}-\mathbf{A}_{0}^{m+1}{{(\mathbf{A}_{0}^{m+1})}^{+}}\mathbf{h}_{:k+1}^{m+1}.
\end{equation}
Finally we can substitute (\ref{c_wave_k1_A0_123}) into (\ref{c_k_final_compute123}) to deduce (\ref{c_wave_ c_k _equal}).
%compare (\ref{c_wave_k1_A0_123}) with (\ref{c_k_final_compute123}) to verify (\ref{c_wave_ c_k _equal}).

Since the condition $\mathbf{C}=\mathbf{0}$ in (\ref{B_Matrix_def1b}) is equivalent to
\begin{equation}\label{c_ 1_p_all0_345}
\mathbf{c}_{:1}^{{}}=\mathbf{c}_{:2}^{{}}=\cdots =\mathbf{c}_{:p}^{{}}=\mathbf{0},
\end{equation}
  we can deduce $\mathbf{\tilde{c}}_{:1}^{{}}=\mathbf{c}_{:1}^{{}}=\mathbf{0}$ from (\ref{c1_wave_c1_equal}), and then we can apply Theorem 2 and (\ref{c_ 1_p_all0_345}) iteratively to deduce $\mathbf{\tilde{c}}_{:2}^{{}}=\mathbf{c}_{:2}^{{}}=\mathbf{0}$, $\mathbf{\tilde{c}}_{:3}^{{}}=\mathbf{c}_{:3}^{{}}=\mathbf{0}$, $\cdots$, and $\mathbf{\tilde{c}}_{:p}^{{}}=\mathbf{c}_{:p}^{{}}=\mathbf{0}$. Correspondingly we have

\underline{\textbf{Theorem 3}}. If $\mathbf{C}=\mathbf{0}$, i.e., %all $p$ columns of $\mathbf{C}$ are zeros,
(\ref{c_ 1_p_all0_345}) satisfies,
then
\begin{equation}\label{c_ wave_1_p_all0_3124}
\mathbf{\tilde{c}}_{:1}^{{}}=\mathbf{\tilde{c}}_{:2}^{{}}=\cdots =\mathbf{\tilde{c}}_{:p}^{{}}=\mathbf{0}.
\end{equation}

%\begin{equation}\label{b_col_def1ab}
%\mathbf{\tilde{b}}_{:k}^{T}=\left\{ \begin{matrix}
%   \mathbf{\tilde{c}}_{:k}^{+} & i{{f}^{{}}}\mathbf{\tilde{c}}_{:k}^{{}}\ne \mathbf{0},  \\
%   {{(1+\mathbf{\tilde{d}}_{:k}^{T}{{{\mathbf{\tilde{d}}}}_{:k}})}^{-1}}\mathbf{\tilde{d}}_{:k}^{T}{{(\mathbf{A}_{k-1}^{m+1})}^{+}} & i{{f}^{{}}}\mathbf{\tilde{c}}_{:k}^{{}}=\mathbf{0}.   \\
%\end{matrix} \right.
%\end{equation}

%$\mathbf{\tilde{b}}_{:k}^{T}=$
%\begin{subnumcases}{\label{b_col_def1ab}}
% \mathbf{\tilde{c}}_{:k}^{+} & i{{f}^{{}}}\mathbf{\tilde{c}}_{:k}^{{}}\ne \mathbf
%{0} &  \label{b_col_def1a}\\
%{{(1+\mathbf{\tilde{d}}_{:k}^{T}{{{\mathbf{\tilde{d}}}}_{:k}})}^{-1}}\mathbf{\tilde{d}}_{:k}^{T}{{(\mathbf{A}_{k-1}^{m+1})}^{+}} & i{{f}^{{}}}\mathbf{\tilde{c}}_{:k}^{{}}=\mathbf{0}. & \label{b_col_def1b}
%\end{subnumcases}

Since ${{\mathbf{\tilde{c}}}_{:k}}=0$ is equivalent to (\ref{RAkm1equ2RAkm1}),
%$\Re (\mathbf{A}_{k}^{m+1})=\Re (\mathbf{A}_{k-1}^{m+1})$ \cite{cite_general_inv_book},
Theorem 3 is also equivalent to:  if $\mathbf{C}=\mathbf{0}$, then
\begin{equation}\label{C_matrix_zero_condi1}
\Re (\mathbf{A}_{p}^{m+1})=\Re (\mathbf{A}_{p-1}^{m+1})=\cdots =\Re (\mathbf{A}_{0}^{m+1}).
\end{equation}
%Theorem 3 is also equivalent to the following Theorem 4.
%
%Theorem 4. If $\mathbf{C}=\mathbf{0}$, i.e.,
%%all $p$ columns of $\mathbf{C}$ are zeros,
%(\ref{c_ 1_p_all0_345}) satisfies,
%then
%%\begin{equation}\label{C_matrix_zero_condi1}
%%\Re (\mathbf{A}_{p}^{m+1})=\Re (\mathbf{A}_{p-1}^{m+1})=\cdots =\Re (\mathbf{A}_{1}^{m+1})=\Re (\mathbf{A}_{0}^{m+1})
%%\end{equation}
%\begin{equation}\label{C_matrix_zero_condi1}
%\Re (\mathbf{A}_{p}^{m+1})=\Re (\mathbf{A}_{p-1}^{m+1})=\cdots =\Re (\mathbf{A}_{0}^{m+1}),
%\end{equation}
%%(where $\mathbf{A}_{0}^{m+1}=\mathbf{A}_{{}}^{m}$),
%i.e., $\Re (\mathbf{A}_{k}^{m+1})=\Re (\mathbf{A}_{k-1}^{m+1})$ always satisfies for $k=1,2,\cdots,p$.

\subsection{Modified Greville's method for BLS Considering All Possible Cases and Improved in Efficiency and Numerical Stability}

From Theorem 3, it can be seen that when the condition in (\ref{B_Matrix_def1b})  satisfies, the condition in (\ref{b_col_def1b}) also satisfies for $k=1,2, \cdots ,p$. Then (\ref{b_col_def1b})  can be applied to compute $\mathbf{\tilde{b}}_{:1}^{T}$, $\mathbf{\tilde{b}}_{:2}^{T}$, $\cdots$, and $\mathbf{\tilde{b}}_{:p}^{T}$. Correspondingly (\ref{b_col_def1b}) can be extended to (\ref{B_Matrix_def1b}).

To improve the numerical stability and reduce the computational complexity,  substitute
(\ref{D_matrix_def_111})
into
(\ref{B_Matrix_def1b})
to obtain
\begin{equation}\label{Bt2DDtAHDtA04325}
{{\bf{B}}^T} ={{({\bf{I}} + {{\bf{D}}^T}{{({{\mathbf{A}}^{m}})}^{+}}{{\mathbf{H}}_{m+1}})}^{ - 1}}{{\bf{D}}^T}{{({{\bf{A}}^m})}^ + },
\end{equation}
which can be written as
%Write (\ref{Bt2DDtAHDtA04325}) as
\begin{equation}\label{DaddDHinvD492573}
{{\bf{B}}^T} ={{({\bf{I}} + {\bf{\tilde D}}{{\mathbf{H}}_{m+1}})}^{ - 1}}{\bf{\tilde D}}
\end{equation}
where
%${\bf{\tilde D}}={{\bf{D}}^T}{{({{\bf{A}}^m})}^ + }$, i.e.,
\begin{equation}\label{DTilde2DtAmInv94832}
{\bf{\tilde D}}={{\bf{D}}^T}{{({{\bf{A}}^m})}^ + }.
\end{equation}
%\begin{equation}\label{DTilde2DtAmInv94832}
%{\bf{\tilde D}}={{\mathbf{H}}_{m+1}^T}\left( {{({{\mathbf{A}}^{m}})}^{+}}\right)^T{{({{\bf{A}}^m})}^ + }.
%\end{equation}
Then we can utilize  equation (20) in \cite{InverseSumofMatrix8312}, i.e.,
 \begin{equation}\label{}
(\mathbf{I}+\mathbf{P} \mathbf{Q})^{-1} \mathbf{P}=\mathbf{P}(\mathbf{I}+\mathbf{Q} \mathbf{P})^{-1},
 \end{equation}
 to deduce
\begin{equation}\label{DaddDHinvD492573Other}
{{\bf{B}}^T} ={\bf{\tilde D}}{{({\bf{I}} + {{\mathbf{H}}_{m+1}}{\bf{\tilde D}})}^{ - 1}}
\end{equation}
from (\ref{DaddDHinvD492573}).
%Since
%%$\mathbf{H}_{m+1}^{{}}$ is $m \times p$
%%%while  from  (\ref{D_matrix_def_111})
%%%and
%%%(\ref{DTilde2DtAmInv94832}) it can be seen that
%%and ${\bf{\tilde D}}$ is $p \times m$.
%$\mathbf{H}_{m+1}^{{}}$
%and ${\bf{\tilde D}}$ are
%$m \times p$ and
% $p \times m$, respectively,
%For the $m \times p$
%$\mathbf{H}_{m+1}^{{}}$ and
%$p \times m$ ${\bf{\tilde D}}$,
%  %Thus
%obviously  (\ref{DaddDHinvD492573Other}) with an $m \times m$ matrix inverse   is more  stable  and efficient
%than (\ref{DaddDHinvD492573}) or (\ref{B_Matrix_def1b})  with a $p \times p$
%matrix
% inverse
%when $m < p$, and
%(\ref{DaddDHinvD492573}) or (\ref{B_Matrix_def1b})  with a $p \times p$
%matrix
% inverse is more stable  and efficient than  (\ref{DaddDHinvD492573Other}) with an $m \times m$ matrix inverse
%  when $m > p$.

$\mathbf{H}_{m+1}^{{}}$ and
${\bf{\tilde D}}$ are $m \times p$  and $p \times m$, respectively.
Then it can be seen that when $m < p$,  (\ref{DaddDHinvD492573Other}) with an $m \times m$ matrix inverse   is more  stable  and efficient
than (\ref{DaddDHinvD492573}) or (\ref{B_Matrix_def1b})  with a $p \times p$
matrix
 inverse.
On the other hand,  when $m > p$,
(\ref{DaddDHinvD492573}) or (\ref{B_Matrix_def1b})  with a $p \times p$
matrix inverse is more stable  and efficient than  (\ref{DaddDHinvD492573Other}) with an $m \times m$ matrix inverse.
Moreover,   when $m>n$,  the computational complexity of ${{\bf{D}}^T}{\bf{D}}$ in (\ref{B_Matrix_def1b}) is lower
 than that of ${\bf{\tilde D}}{{\mathbf{H}}_{m+1}}$  in  (\ref{DaddDHinvD492573}).

On the other hand, let us consider the condition in (\ref{B_Matrix_def1a}), i.e.,
$\mathbf{C}\ne 0$. Firstly, let us give

\underline{\textbf{The Inverse Negative Proposition of Theorem 3}}.   If
(\ref{c_ wave_1_p_all0_3124})
is not satisfied,
i.e.,
there is at least one $\mathbf{\tilde{c}}_{:k}^{{}}\ne \mathbf{0}$ ($1\le k\le p$), then $\mathbf{C}\ne 0$.
%which is equivalent to  $\Re (\mathbf{A}_{k}^{m+1})\ne \Re (\mathbf{A}_{k-1}^{m+1})$.
%for at least one $k$.
%When $\mathbf{C}\ne 0$, from Theorems 3 and 4, it can be seen that there is at least one $\mathbf{\tilde{c}}_{:k}^{{}}\ne \mathbf{0}$, i.e., $\Re (\mathbf{A}_{k}^{m+1})\ne \Re (\mathbf{A}_{k-1}^{m+1})$ for at least one $k$, where $1\le k\le p$.

To compute ${{\mathbf{B}}^{T}}$ by (\ref{B_Matrix_def1a}), obviously the condition in (\ref{b_col_def1a})
 (i.e., $\mathbf{\tilde{c}}_{:k}^{{}}\ne \mathbf{0}$)  should be satisfied for all $k=1,2,\cdots,p$.
%, for all $k=1,2,\cdots,p$.
%(\ref{B_Matrix_def1a}) ¸Ä³ÉÕýÈ·µÄ¸üÇ¿µÄÌõ¼þ£¬(\ref{B_Matrix_def1b})µÄÌõ¼þûÓÐ´í£¬´ËÍ⣬ÕâÁ½¸öÌõ¼þûÓаüÀ¨ËùÓеÄÌõ¼þ£¬ËùÒÔ±ØÐë°ÑÕâÁ½¸öÌõ¼þÖ®ÍâÆäËü¿ÉÄܵÄÌõ¼þÒ²¿¼ÂǽøÀ´£¬i.e.,
%(\ref{B_Matrix_def1a}) ¸Ä³ÉÕýÈ·µÄ¸üÇ¿µÄÌõ¼þ£¬
This condition is much stronger than the above-described condition of at least one $\mathbf{\tilde{c}}_{:k}^{{}}\ne \mathbf{0}$. Thus in (\ref{B_Matrix_def1a}), ``$if \quad \mathbf{C}\ne 0$'' should be modified into
%¡°if all p $\mathbf{\tilde{c}}_{:k}^{{}}$s satisfy $\mathbf{\tilde{c}}_{:k}^{{}}\ne \mathbf{0}$ ($1\le k\le p$)¡±,
``if each $\mathbf{\tilde{c}}_{:k}^{{}}\ne \mathbf{0}$ ($1\le k\le p$)'',
%(\ref{B_Matrix_def1b})µÄÌõ¼þûÓÐ´í£¬´ËÍ⣬ÕâÁ½¸öÌõ¼þûÓаüÀ¨ËùÓеÄÌõ¼þ£¬ËùÒÔ±ØÐë°ÑÕâÁ½¸öÌõ¼þÖ®ÍâÆäËü¿ÉÄܵÄÌõ¼þÒ²¿¼ÂǽøÀ´£¬
and it is required to consider the
%cases that only
condition of only
%part of the $p$ $\mathbf{\tilde{c}}_{:k}^{{}}$s satisfy $\mathbf{\tilde{c}}_{:k}^{{}}\ne \mathbf{0}$ ($1\le k\le p$), i.e.,
$i$ ($1\le i\le p-1$) $\mathbf{\tilde{c}}_{:k}$s satisfying $\mathbf{\tilde{c}}_{:k}=\mathbf{0}$, i.e., ${\bf{C}} \ne {\bf{0}}$  but several ${\bf{\tilde c}}_{:k} = {\bf{0}}$.

Now we can modify
%Correspondingly
(\ref{B_Matrix_def1ab})
%should be modified
into
\begin{subequations}{\label{B_Matrix_def2abc}}
 \begin{numcases}
{  {{\bf{B}}^T} = }
  {{({\bf{I}} + {{\bf{D}}^T}{\bf{D}} )^{ - 1}}}
{\bf{\tilde D}} \ if \,  {\bf{C}} = {\bf{0}}, m \ge max(n,p)     &  \label{B_Matrix_def2bOld} \\
{{({\bf{I}} + {\bf{\tilde D}}{{\mathbf{H}}_{m+1}})}^{ - 1}}{\bf{\tilde D}} \ if\ {\bf{C}} = {\bf{0}}, n \ge m \ge p  &  \label{B_Matrix_def2bOther} \\
{\bf{\tilde D}}{{({\bf{I}} + {{\mathbf{H}}_{m+1}}{\bf{\tilde D}})}^{ - 1}} \quad \quad \ if\ {\bf{C}} = {\bf{0}},  m \le p  &  \label{B_Matrix_def2b} \\
 {{{\bf{C}}^ + }\quad \quad \quad \quad \ if\ each\ {\bf{\tilde c}}_{:k}^{} \ne {\bf{0}}(1 \le k \le p)}  &  \label{B_Matrix_def2a} \\
{\cdots \quad \quad \quad \  if\ {\bf{C}} \ne {\bf{0}} \ but\ several\ {\bf{\tilde c}}_{:k}^{} = {\bf{0}}},  &  \label{B_Matrix_def2c}
\end{numcases}
\end{subequations}
where the original Greville's method \cite{cite_general_inv_book} can be utilized
to compute ${{\left( \mathbf{A}_{p}^{m+1} \right)}^{+ }}={{\left( \mathbf{A}_{{}}^{m+1} \right)}^{+ }}$ from ${{\left( \mathbf{A}_{0}^{m+1} \right)}^{+ }}={{\left( \mathbf{A}_{{}}^{m} \right)}^{+ }}$ by $p$ iterations,
if the condition for (\ref{B_Matrix_def2c}) is satisfied.
%${\bf{\tilde D}}$ is computed by (\ref{D_matrix_def_111}) and (\ref{DTilde2DtAmInv94832}),
%(\ref{DaddDHinvD492573Other})  and (\ref{DaddDHinvD492573}) are included,
% tempBetaD= beta_paper*Dt_paper';
%t2notManySamples=1;
%if size(t2_func,1)<=size(t2_func,2)
%b_func=tempBetaD*inv(eye(size(d_t_func,1))+t2_func*tempBetaD);
%% B_paper=tempBetaD*inv(eye(size(Dt_paper,1))+t2*tempBetaD);
%else
%b_func=inv(eye(size(beta,1))+ tempBetaD*t2_func)*tempBetaD;
%% B_paper=inv(eye(size(beta_paper,1))+ tempBetaD*t2)*tempBetaD;
%end
%\begin{subnumcases}{\label{B_Matrix_def2abc}}
%{{{\bf{B}}^T} = {{\bf{C}}^ + }\quad if\ each\ {\bf{\tilde c}}_{:k}^{} \ne {\bf{0}}(1 \le k \le p)} &  \label{B_Matrix_def2a}\\
%{{{\bf{B}}^T} = {{({\bf{I}} + {{\bf{D}}^T}{\bf{D}})}^{ - 1}}{{\bf{D}}^T}{{({{\bf{A}}^m})}^ + }\quad if\ {\bf{C}} = {\bf{0}}}  &  \label{B_Matrix_def2b}\\
%{{{\bf{B}}^T} = \cdots \quad if\ {\bf{C}} \ne {\bf{0}} \ but\ several\ {\bf{\tilde c}}_{:k}^{} = {\bf{0}}}, & \label{B_Matrix_def2c}
%\end{subnumcases}

%In (B_Matrix_def2c)£¬»Øµ½Ò»ÁÐÒ»ÁеÄË㣻
%»òÕߣ¬¸ü¸ßЧµÄËã·¨ÊÇ£¬ÕÒµ½µÚÒ»¸öΪ0µÄcÖ®ºó£¬
%1.Èç¹û´ËÇ°Óв»ÊÇÁãµÄ¶à¸öÁУ¬¶ÔËüÃÇÓÃa£¬Ëã³ö¹ãÒåÄæ¾ØÕó
%2.¶ÔÓÚµ±Ç°Îª0µÄÁУ¬
%2.1Èç¹ûC¿ÉÒÔÖظ´ÀûÓ㨼´µ±Ç°Îª0µÄÁÐÊǵ±Ç°CµÄµÚ1ÁУ©£¬¿´¿´CºóÃ漸ÁÐÓÐûÓнô°¤×ŵÄ0£¬ÓУ¬Ôòµ±Ç°ÁÐÓëºóÃ漸ÁÐÒ»Æð£¬ÓÃb Ëã¹ãÒåÄæ¾ØÕó£»Ã»ÓУ¬ÔòÈçͬ2.2£¬Óõ¥ÁеÄËã·¨Çó¹ãÒåÄæ¾ØÕó
%2.2Èç¹ûC²»¿ÉÒÔÖظ´ÀûÓ㨼´µ±Ç°Îª0µÄÁв»Êǵ±Ç°CµÄµÚ1ÁУ©£¬Óõ¥ÁеÄËã·¨Çó¹ãÒåÄæ¾ØÕó
%¶ÔÓÚHÓàϵÄÁУ¬¼ÌÐøÒÔÉϵÄËã·¨£¬Ö±µ½HËùÓеÄÁУ¬¶ÔÓ¦µÄc_waveÅжÏÍê±Ï¡£
%
%Õû¸öËã·¨£º°Ñ´úÂë·­Òë¹ýÀ´¼´¿É
%
%wwwwwwwwwwwwwwwwwwwwwwwwwwwww

\subsection{To Apply the Recently Proposed Inverse Cholesky Factorization to Compute All ${{\mathbf{\tilde{c}}}_{:k}}$s Efficiently}
% Modified Greville's method for BLS Considering All Possible Cases and Improved in Efficiency and Numerical Stability

In
(\ref{B_Matrix_def2a})
and
(\ref{B_Matrix_def2c}),
 all $p$ ${{\mathbf{\tilde{c}}}_{:k}}$s ($k=1,2,\cdots,p$) are required.  If they are computed by
%(\ref{c_col_define1})
(\ref{c_wave_k_final_compute123}) in $p$ iterations, $p-1$ ${{(\mathbf{A}_{k-1}^{m+1})}^{+}}$s ($k=2,3,\cdots,p$) also need to be computed by
 (\ref{A_1col_inv_book1})
% , (\ref{d_col_define4311}) and
%(\ref{b_col_def1ab})
 in $p-1$ iterations, and then actually it is no longer required to apply (\ref{A_1col_inv_book1matrix}),
 %() and ()
 (\ref{D_matrix_def_111}),
(\ref{C_matrix_def_111}) and
(\ref{B_Matrix_def1ab})
 once to compute ${{(\mathbf{A}_{p}^{m+1})}^{+}}$ from ${{(\mathbf{A}_{0}^{m+1})}^{+}}$ directly. Thus in what follows, we will propose an efficient algorithm to compute all $p$ ${{\mathbf{\tilde{c}}}_{:k}}$s efficiently, which is based on the recently proposed efficient inverse Cholesky factorization~\cite{my_inv_chol_paper}, and does not require the above-mentioned $p-1$ ${{(\mathbf{A}_{k-1}^{m+1})}^{+}}$s ($k=2,3,\cdots,p$).

To apply  the efficient inverse Cholesky factorization~\cite{my_inv_chol_paper}, firstly let us derive

\underline{\textbf{Theorem 4}}. If each ${\bf{\tilde c}}_{:k}^{} \ne {\bf{0}}$ ($1 \le k \le p$), then
%$\mathbf{C}$ is full column rank, and
 $\mathbf{C}_{{}}^{T}\mathbf{C}$ is positive definite.
%, and vice versa.
%4. ËãC pseudoinverse, ÓÃinv(C¡¯C)¡£

Proof.   If each ${\bf{\tilde c}}_{:k}^{} \ne {\bf{0}}(1 \le k \le p)$, the condition in (\ref{B_Matrix_def2a}) is satisfied,
%i.e., all $p$ $\mathbf{\tilde{c}}_{:k}$s satisfy $\mathbf{\tilde{c}}_{:k}^{{}}\ne {{\mathbf{0}}^{{}}}(1\le k\le p)$.
%   we have
%\begin{equation}
%\label{B_C_PR_inv123}{{\mathbf{B}}^{T}}={{\mathbf{C}}^{+}}.
%\end{equation}
%From (\ref{C_matrix_def_111}), i.e., $\mathbf{C}={{\mathbf{H}}_{m+1}}-{{\mathbf{A}}^{m}}\mathbf{D}$, it can be seen that the row size is usually bigger than the column size
%%???
%and then ${{\mathbf{C}}^{+}}$ can be written as the left inverse, i.e., (\ref{B_C_PR_inv123})
and then
%(\ref{B_Matrix_def2a}) can be utilized,
we can utilize the computation in (\ref{B_Matrix_def2a}),
 which
can be written as~\cite{BL_trans_paper}
\begin{equation}\label{C_left_inv_def1}
{{\mathbf{B}}^{T}}={{\mathbf{C}}^{+}}={{({{\mathbf{C}}^{T}}\mathbf{C})}^{-1}}{{\mathbf{C}}^{T}}.
\end{equation}
%5. inv C by inv chol, ZFµÄÍƵ¼¡£
%We can apply the efficient Cholesky factorization recently proposed in \cite{my_inv_chol_paper} to compute the inverse Cholesky factor of ${{\mathbf{C}}^{T}}\mathbf{C}$ in (\ref{C_left_inv_def1}), i.e., the upper-triangular matrix $\mathbf{G}$ that satisfies
%$\mathbf{FF}_{{}}^{T}={{(\mathbf{C}_{{}}^{T}\mathbf{C})}^{-1}}$.
Since ${{(\mathbf{C}_{{}}^{T}\mathbf{C})}^{-1}}$ exists,  $\mathbf{C}$ must be full column rank, and then $\mathbf{C}_{{}}^{T}\mathbf{C}$ must be positive definite~\cite{Matrix_Computations_book}.
% On the other hand, if $\mathbf{C}_{{}}^{T}\mathbf{C}$ is positive definite, there exists inverse Cholesky of $\mathbf{C}_{{}}^{T}\mathbf{C}$, and then cWave computed by the following () not zero.

From Theorem 4, it can be seen that if each ${\bf{\tilde c}}_{:k}^{} \ne {\bf{0}}(1 \le k \le p)$,
%$\mathbf{C}_{{}}^{T}\mathbf{C}$ is positive definite
%Thus
there exists~\cite{Matrix_Computations_book}
the  Cholesky factor of the positive definite  $\mathbf{C}_{{}}^{T}\mathbf{C}$, i.e.,  the lower-triangular
${\bf{\Omega }}$ that satisfies
\begin{equation}\label{OmegaCholR89866}
{\bf{\Omega}}{\bf{\Omega }}^T=\mathbf{C}_{{}}^{T}\mathbf{C},
 \end{equation}%[]
from
%(\ref{OmegaCholR89866})
which we can deduce
\begin{equation}\label{OmegaCholR89866Inv}
{\bf{\Omega }}^{-T}{{\bf{\Omega }}^{-1}}={{(\mathbf{C}^{T}{{\mathbf{C}}})}^{-1}}.
\end{equation}
%and then
From (\ref{OmegaCholR89866Inv}) it can be seen that
the upper-triangular  ${\bf{\Omega }}^{-T}$ is
the inverse Cholesky factor~\cite{my_inv_chol_paper} of $\mathbf{C}_{{}}^{T}\mathbf{C}$,
%i.e.,
 which can be denoted as
\begin{equation}\label{G2OmegaT9498}
{\mathbf{G}}={\bf{\Omega }}^{-T}.
 \end{equation}%[]
% To compute ${{(\mathbf{C}_{{}}^{T}\mathbf{C})}^{-1}}$ in (\ref{C_left_inv_def1}) efficiently,
% %Let. Then the inverse Cholesky factor of $\mathbf{C}_{k}^{T}{{\mathbf{C}}_{k}}$ is
% we can utilize ${{\mathbf{G}}_{k}}$ that satisfies
%\begin{equation}\label{R_def12445}
%{{\mathbf{G}}_{k}}\mathbf{G}_{k}^{T}={{(\mathbf{C}_{k}^{T}{{\mathbf{C}}_{k}})}^{-1}},
%\end{equation}
%where ${{\mathbf{C}}_{k}}$ denote the first $k$ columns of ${\mathbf{C}}$, and the upper-triangular ${{\mathbf{G}}_{k}}$ is the inverse Cholesky factor~\cite{my_inv_chol_paper} of $\mathbf{C}_{k}^{T}{{\mathbf{C}}_{k}}$.

To obtain ${\mathbf{G}}$, we can utilize the efficient Cholesky factorization proposed in \cite{my_inv_chol_paper} to compute ${{\mathbf{G}}_{k}}$ from ${{\mathbf{G}}_{k-1}}$ iteratively for $k=2,3,\cdots,p$, by equation (11) in \cite{my_inv_chol_paper}, i.e.,
\begin{equation}\label{F_iterative_Add}
%{{\mathbf{G}}_{k}}=\left[ \begin{matrix}
%   {{\mathbf{G}}_{k-1}} & {{\mathbf{u}}_{k-1}}  \\
%   \mathbf{0}_{k-1}^{T} & {{\eta }_{k}}  \\
%\end{matrix} \right],
{{\mathbf{G}}_{k}}=\left[ \begin{matrix}
   {{\mathbf{G}}_{k-1}} & {{{\mathbf{\tilde{g}}}}_{:k}}  \\
   \mathbf{0}_{k-1}^{T} & {{g}_{kk}}  \\
\end{matrix} \right],
\end{equation}
where
\begin{equation}\label{ita_myChol_1212}
{{g}_{kk}}=1/\sqrt{\mathbf{c}_{:k}^{T}{{\mathbf{c}}_{:k}}-\mathbf{c}_{:k}^{T}\mathbf{C}_{k-1}^{{}}{{\mathbf{G}}_{k-1}}\mathbf{G}_{k-1}^{T}\mathbf{C}_{k-1}^{T}{{\mathbf{c}}_{:k}}}
\end{equation}
% {{\eta }_{k}}
and
\begin{equation}\label{u_myChol_12321}
%{{\mathbf{u}}_{k-1}}=-{{\eta }_{k}}{{\mathbf{G}}_{k-1}}\mathbf{G}_{k-1}^{T}\mathbf{C}_{k-1}^{T}{{\mathbf{c}}_{:k}}.
{{{\mathbf{\tilde{g}}}}_{:k}}=-{{g}_{kk}}{{\mathbf{G}}_{k-1}}\mathbf{G}_{k-1}^{T}\mathbf{C}_{k-1}^{T}{{\mathbf{c}}_{:k}}.
\end{equation}
%\begin{subnumcases}{\label{ita_myChol_u_myChol}}
% {{\eta }_{k}}=1/\sqrt{\mathbf{c}_{:k}^{T}{{\mathbf{c}}_{:k}}-\mathbf{c}_{:k}^{T}\mathbf{C}_{k-1}^{{}}{{\mathbf{G}}_{k-1}}\mathbf{G}_{k-1}^{T}\mathbf{C}_{k-1}^{T}{{\mathbf{c}}_{:k}}}&  \label{ita_myChol_1212}\\
%{{\mathbf{u}}_{k-1}}=-{{\eta }_{k}}{{\mathbf{G}}_{k-1}}\mathbf{G}_{k-1}^{T}\mathbf{C}_{k-1}^{T}{{\mathbf{c}}_{:k}}. & \label{u_myChol_12321}
%\end{subnumcases}
%\begin{equation}\label{ita_myChol_1212}
%
%\end{equation}(ita_myChol_1212) ,
%\begin{equation}\label{u_myChol_12321}
%
%\end{equation}(u_myChol_12321).
%(The above 2 should be in a group).
Notice that equations (\ref{ita_myChol_1212}) and (\ref{u_myChol_12321})
%Equation (\ref{ita_myChol_u_myChol})
 %can be
 are derived from  equations (3) and (17) in \cite{my_inv_chol_paper}.

Now let us consider the case that
 ${{\mathbf{\tilde{c}}}_{:i}}\ne 0$ is satisfied for only
$i=1,2,\cdots, k$, where $k<p$. According to Theorem 4, we can conclude that
 $\mathbf{C}_{{k}}^{T}\mathbf{C}_{{k}}$ is positive definite, where ${{\mathbf{C}}_{k}}$
 denotes the first $k$ columns of ${\mathbf{C}}$.
%Accordingly,  denote  as.
With the positive definite $\mathbf{C}_{{k}}^{T}\mathbf{C}_{{k}}$,
 %, and (\ref{R_def12445}).
%Obviously,
 the upper-triangular ${{\mathbf{G}}_{k}}$ in  (\ref{F_iterative_Add}) can be computed by~\cite{my_inv_chol_paper}
\begin{equation}\label{R_def12445}
{{\mathbf{G}}_{k}}\mathbf{G}_{k}^{T}={{(\mathbf{C}_{k}^{T}{{\mathbf{C}}_{k}})}^{-1}},
\end{equation}
and
 $\mathbf{B}_{k}^{T}$ %corresponding to ${{\mathbf{C}}_{k}}$
 can be computed by
%applying
% in Theorem 4,
%can be applied to compute
$\mathbf{B}_{k}^{T}=\mathbf{C}_{k}^{+ }={{({{\mathbf{C}}_{k}^{T}}\mathbf{C}_{k})}^{-1}}{{\mathbf{C}}_{k}^{T}}$ (i.e., (\ref{C_left_inv_def1})),
into which substitute (\ref{R_def12445}) to obtain
\begin{equation}\label{comments_myB_def111}
\mathbf{B}_{k}^{T}=\mathbf{C}_{k}^{+ }={{\mathbf{G}}_{k}}\mathbf{G}_{k}^{T}\mathbf{C}_{k}^{T}.
\end{equation}
In (\ref{comments_myB_def111}), $\mathbf{C}_{k}$ can be computed by (\ref{C_matrix_def_111}), i.e.,
\begin{equation}\label{C_matrix_def_onlypartH}
{{\mathbf{C}}_{k}}=\mathbf{H}_{m+1}^{:,1:k}-{{\mathbf{A}}^{m}}{{\mathbf{D}}_{k}},
\end{equation}
where  ${{\mathbf{D}}_{k}}$ is computed by (\ref{D_matrix_def_111}), i.e.,
\begin{equation}\label{DjMinus1EquAH}
{{\mathbf{D}}_{k}}= {{({{\mathbf{A}}^{m}})}^{+}}\mathbf{H}_{m+1}^{:,1:k}\text{ }.
\end{equation}
%Finally,
 $\mathbf{B}_{k}$ and  ${{\mathbf{D}}_{k}}$ can be applied to
% update  the
% pseudoinverse
 compute ${{(\mathbf{A}_{k}^{m+1})}^{+}}$
  by  (\ref{A_1col_inv_book1matrix}), i.e.,
%Finally, $\mathbf{B}_{k}$ and  ${{\mathbf{D}}_{k}}$ can be applied to update  the
% pseudoinverse by  (\ref{A_1col_inv_book1matrix}), i.e.,
%\begin{multline}\label{BLK_general_inv_def1223}
%{{(\mathbf{A}_{k}^{m+1})}^{+}}={{\left[ {{\mathbf{A}}^{m}}|\mathbf{H}_{m+1}^{:,1:k} \right]}^{+}} \\
%=\left[ \begin{matrix}
%   {{({{\mathbf{A}}^{m}})}^{+}}-{{\mathbf{D}}_{k}}\mathbf{B}_{k}^{T}  \\
%   \mathbf{B}_{k}^{T}  \\
%\end{matrix} \right]\text{ }.
%\end{multline}
\begin{equation}\label{BLK_general_inv_def1223}
{{(\mathbf{A}_{k}^{m+1})}^{+}}={{\left[ {{\mathbf{A}}^{m}}|\mathbf{H}_{m+1}^{:,1:k} \right]}^{+}} =\left[ \begin{matrix}
   {{({{\mathbf{A}}^{m}})}^{+}}-{{\mathbf{D}}_{k}}\mathbf{B}_{k}^{T}  \\
   \mathbf{B}_{k}^{T}  \\
\end{matrix} \right]\text{ }.
\end{equation}

In the above case of
 ${{\mathbf{\tilde{c}}}_{:i}}\ne 0$
 %is satisfied
  for
  %only
$i=1,2,\cdots, k-1$,
%where $k<p$.
the corresponding $\mathbf{C}_{k-1}$ and ${{\mathbf{G}}_{k-1}}$
can be utilized to compute ${{\mathbf{\tilde{c}}}_{:k}}$ and its squared length
efficiently, as shown in the following Theorem 5. The proof of Theorem 5 is given
% from the above
%-mentioned
  in Appendix A. % we will prove

\underline{\textbf{Theorem 5}}. When ${{\mathbf{\tilde{c}}}_{:i}}\ne 0$ for all $1\le i\le k-1$, ${{\mathbf{\tilde{c}}}_{:k}}$ ($k=2,3,\cdots, p$) defined by (\ref{c_wave_k_final_compute123}) is equal to
\begin{equation}\label{c_col_define4myChol}
{{\mathbf{\tilde{c}}}_{:k}}={{\mathbf{c}}_{:k}}-\mathbf{C}_{k-1}^{{}}{{\mathbf{G}}_{k-1}}\mathbf{G}_{k-1}^{T}\mathbf{C}_{k-1}^{T}{{\mathbf{c}}_{:k}},
\end{equation}
and the squared length of ${{\mathbf{\tilde{c}}}_{:k}}$  can be computed by
\begin{multline}\label{c_lentghcCFFcc}
\left|{{\mathbf{\tilde{c}}}_{:{k}}} \right|^2=\mathbf{\tilde{c}}_{:k}^{T}{{\mathbf{\tilde{c}}}_{:k}} \\
=\mathbf{c}_{:k}^{T}{{\mathbf{c}}_{:k}}-\mathbf{c}_{:k}^{T}\mathbf{C}_{k-1}^{{}}
{{\mathbf{G}}_{k-1}}\mathbf{G}_{k-1}^{T}\mathbf{C}_{k-1}^{T}{{\mathbf{c}}_{:k}}.
\end{multline}

It can be seen that  the pseudoinverse
 ${{(\mathbf{A}_{k-1}^{m+1})}^{+}}$  utilized in (\ref{c_wave_k_final_compute123}) to compute $\mathbf{\tilde{c}}_{:k}^{{}}$
 is no longer required in (\ref{c_col_define4myChol}) and (\ref{c_lentghcCFFcc}).
  % and it is possible to compute
%  $\mathbf{A}_{p}^{m+1}$ from $\mathbf{A}_{0}^{m+1}$ by just one iteration.
  Moreover, we can substitute (\ref{c_lentghcCFFcc})
   into (\ref{ita_myChol_1212}) to compute ${{g}_{kk}}$ by
\begin{equation}\label{ita_myChol4c_wave}
{{g}_{kk}}=1/\sqrt{\left|  {{\mathbf{\tilde{c}}}_{:{k}}} \right|^2}= 1/\sqrt{\mathbf{\tilde{c}}_{:k}^{T}{{\mathbf{\tilde{c}}}_{:k}}}.
\end{equation}

With Theorem 5,
%Moreover,
we can %also
 prove the Inverse Proposition of Theorem 4, i.e.,

\underline{\textbf{Theorem 6}}. If $\mathbf{C}_{{}}^{T}\mathbf{C}$ is positive definite, then each ${{\mathbf{\tilde{c}}}_{:k}}$ ($1 \le k \le p$) defined by
(\ref{c_wave_k_final_compute123})
 %(\ref{c_col_define1})
  satisfies
${\bf{\tilde c}}_{:k}^{} \ne {\bf{0}}$.
%, and vice versa.
%4. ËãC pseudoinverse, ÓÃinv(C¡¯C)¡£

Proof.   If $\mathbf{C}_{{}}^{T}\mathbf{C}$ is positive definite, there exists~\cite{Matrix_Computations_book}
the  Cholesky factor of $\mathbf{C}_{{}}^{T}\mathbf{C}$, i.e.,  ${\bf{\Omega }}$ satisfying (\ref{OmegaCholR89866}),
and then there also exists the inverse Cholesky factor of $\mathbf{C}_{{}}^{T}\mathbf{C}$, i.e., ${\mathbf{G}}$ satisfying (\ref{G2OmegaT9498}).
Accordingly, we can compute
${{\mathbf{G}}_{k}}$ from ${{\mathbf{G}}_{k-1}}$ iteratively for
$k=2, 3, \cdots, p$
 %can be computed
 by
 %(\ref{ita_myChol4c_wave}),
(\ref{ita_myChol4c_wave}), (\ref{u_myChol_12321}) and (\ref{F_iterative_Add}),
and the initial ${{\mathbf{G}}_{1}}$ can be computed by
% From (\ref{R_def12445}), we can deduce that the initial
  \begin{equation}\label{initialF1compute}
{{\mathbf{G}}_{1}}=1/\sqrt{\mathbf{c}_{:1}^{T}{{\mathbf{c}}_{:1}}}
 \end{equation}
 that is deduced from (\ref{R_def12445}).
%(\ref{initialF1compute}),.
From (\ref{initialF1compute}) we obtain ${\bf{c}}_{:1}^{} \ne {\bf{0}}$, from which and (\ref{c1_wave_c1_equal}) we deduce
${\bf{\tilde c}}_{:1}^{} \ne {\bf{0}}$. Moreover, from (\ref{ita_myChol4c_wave}) we deduce
${\bf{\tilde c}}_{:k}^{} \ne {\bf{0}}$ ($2 \le k \le p$) for ${\bf{\tilde c}}_{:k}$ computed by (\ref{c_col_define4myChol}).
%(where $\mathbf{c}_{:1}^{{}} = \mathbf{\tilde{c}}_{:1}^{{}}$)
%Then all the ${\bf{\tilde c}}_{:k}$s in the above-mentioned (\ref{initialF1compute}),  (\ref{ita_myChol4c_wave}), (\ref{u_myChol_12321}) and (\ref{F_iterative_Add})
%must satisfy
%, which can be
% or (for $k=1$).
Since
%${\bf{\tilde c}}_{:1}$ computed by (\ref{c1_wave_c1_equal}) satisfies
 ${\bf{\tilde c}}_{:1}^{} \ne {\bf{0}}$,
from Theorem 5 we can deduce  that ${{\mathbf{\tilde{c}}}_{:2}}$ defined by
(\ref{c_wave_k_final_compute123})
%(\ref{c_col_define1})
 is equal to
${{\mathbf{\tilde{c}}}_{:2}}$ computed by (\ref{c_col_define4myChol}), and then is not zero. Similarly, we can apply Theorem 5 iteratively
to deduce that ${{\mathbf{\tilde{c}}}_{:k}}$ defined by
(\ref{c_wave_k_final_compute123})
%(\ref{c_col_define1})
 is equal to
${{\mathbf{\tilde{c}}}_{:k}}$ computed by (\ref{c_col_define4myChol}), and then is not zero, for $k=3, 4,\cdots, p$.

Now from Theorem 4 and Theorem 6, it can be seen that the condition of $\mathbf{C}_{{}}^{T}\mathbf{C}$ being positive definite is equivalent to
the condition of each ${\bf{\tilde c}}_{:k}^{} \ne {\bf{0}}$  ($1 \le k \le p$), where ${{\mathbf{\tilde{c}}}_{:k}}$ is defined by
(\ref{c_wave_k_final_compute123}).
Then
we can write (\ref{B_Matrix_def2abc}) as
%\begin{subequations}{\label{Bt3Choice98476ColRankPosDef}}
% \begin{numcases}
%{  {{\bf{B}}^T} = }
%  {{({\bf{I}} + {{\bf{D}}^T}{\bf{D}} )^{ - 1}}}
%{\bf{\tilde D}} \ if \,  {\bf{C}} = {\bf{0}}, m \ge max(n,p)     &  \label{B_Matrix_def2bNEWColRankOriginalPosDef} \\
%{{({\bf{I}} + {\bf{\tilde D}}{{\mathbf{H}}_{m+1}})}^{ - 1}}{\bf{\tilde D}} \ if\ {\bf{C}} = {\bf{0}}, n \ge m \ge p  &  \label{B_Matrix_def2bNEWColRankOtherPosDef} \\
%{\bf{\tilde D}}{{({\bf{I}} + {{\mathbf{H}}_{m+1}}{\bf{\tilde D}})}^{ - 1}} \quad \quad \ if\ {\bf{C}} = {\bf{0}},  m \le p  &  \label{B_Matrix_def2bNEWColRankPosDef} \\
%{ {{\bf{C}}^ + }\quad  \quad \quad  \  if\ \mathbf{C}_{{}}^{T}\mathbf{C} \ is \  positive\  definite} &  \label{B_Matrix_def2aNEWColRankPosDef} \\
%{\cdots \quad \ if\ {{\bf{C}} \ne 0}\ and \ \mathbf{C}_{{}}^{T}\mathbf{C} \ is\ not \  positive\  definite}. & \label{B_Matrix_def2cNEWColRankPosDef}
%\end{numcases}
%\end{subequations}

\begin{subequations}{\label{Bt3Choice98476ColRankPosDef}}
 \begin{numcases}
{  {{\bf{B}}^T} = }
  {{({\bf{I}} + {{\bf{D}}^T}{\bf{D}} )^{ - 1}}}
{\bf{\tilde D}} \ if \,  {\bf{C}} = {\bf{0}}, m \ge max(n,p)     &  \label{B_Matrix_def2bNEWColRankOriginalPosDef} \\
{{({\bf{I}} + {\bf{\tilde D}}{{\mathbf{H}}_{m+1}})}^{ - 1}}{\bf{\tilde D}} \ if\ {\bf{C}} = {\bf{0}}, n \ge m \ge p  &  \label{B_Matrix_def2bNEWColRankOtherPosDef} \\
{\bf{\tilde D}}{{({\bf{I}} + {{\mathbf{H}}_{m+1}}{\bf{\tilde D}})}^{ - 1}} \quad \quad \ if\ {\bf{C}} = {\bf{0}},  m \le p  &  \label{B_Matrix_def2bNEWColRankPosDef} \\
{ {{\bf{C}}^ + }\quad  \quad \quad  \  \  if\ \mathbf{C}_{{}}^{T}\mathbf{C} \ is \  positive\  definite} &  \label{B_Matrix_def2aNEWColRankPosDef} \\
\cdots \cdots  \quad \quad \quad \begin{array}{l}
 \ if\ {{\bf{C}} \ne 0}\ and \ \mathbf{C}_{{}}^{T}\mathbf{C} \ is  \\
\ not \  positive\  definite
\end{array} \ . & \label{B_Matrix_def2cNEWColRankPosDef}
\end{numcases}
\end{subequations}

%{}

%(\ref{c_col_define1}).
On the other hand,  the condition of $\mathbf{C}$ being full column rank is equivalent to
 the condition of $\mathbf{C}_{{}}^{T}\mathbf{C}$ being positive definite~\cite{Matrix_Computations_book}, and then is also equivalent to
 the condition of each ${\bf{\tilde c}}_{:k}^{} \ne {\bf{0}}(1 \le k \le p)$.
%if $\mathbf{C}$ is full column rank, $\mathbf{C}_{{}}^{T}\mathbf{C}$ must be positive definite~\cite{Matrix_Computations_book}, and
%then from Theorem 6 we can deduce that each ${{\mathbf{\tilde{c}}}_{:k}}\ne {\bf{0}}$. On the other hand,
%if each ${{\mathbf{\tilde{c}}}_{:k}}\ne {\bf{0}}$, from Theorem 4 we have $\mathbf{C}$ is full column rank. Thus
%the condition of each ${\bf{\tilde c}}_{:k}^{} \ne {\bf{0}}(1 \le k \le p)$, and
%Accordingly,
Then
we can also write (\ref{B_Matrix_def2abc}) as
\begin{subequations}{\label{Bt3Choice98476ColRank}}
 \begin{numcases}
{  {{\bf{B}}^T} = }
  {{({\bf{I}} + {{\bf{D}}^T}{\bf{D}} )^{ - 1}}}
{\bf{\tilde D}} \ if \,  {\bf{C}} = {\bf{0}}, m \ge max(n,p)     &  \label{B_Matrix_def2bNEWColRankOriginal} \\
{{({\bf{I}} + {\bf{\tilde D}}{{\mathbf{H}}_{m+1}})}^{ - 1}}{\bf{\tilde D}} \ if\ {\bf{C}} = {\bf{0}}, n \ge m \ge p  &  \label{B_Matrix_def2bNEWColRankOther} \\
{\bf{\tilde D}}{{({\bf{I}} + {{\mathbf{H}}_{m+1}}{\bf{\tilde D}})}^{ - 1}} \quad \quad \ if\ {\bf{C}} = {\bf{0}},  m \le p  &  \label{B_Matrix_def2bNEWColRank} \\
{ {{\bf{C}}^ + }\quad  \quad \quad \quad \  if\ \mathbf{C} \ is \  full \ column \ rank} &  \label{B_Matrix_def2aNEWColRank} \\
{\cdots \quad if\ {\bf{C}} \ne {\bf{0}} \ is\ not \  full \ column \ rank}. & \label{B_Matrix_def2cNEWColRank}
\end{numcases}
\end{subequations}
%{ {{({\bf{I}} + {{\bf{D}}^T}{\bf{D}})}^{ - 1}}{{\bf{D}}^T}{{({{\bf{A}}^m})}^ + }\quad if\ {\bf{C}} = {\bf{0}}}  &  \label{B_Matrix_def2bNEWColRank}\\

%\section{More and Most are my words}
%
%
%When the condition in (\ref{B_Matrix_def2c})
%%(i.e., that in (\ref{B_Matrix_def2cNEW}))
% satisfies, we can simply go back to the original Greville's method \cite{cite_general_inv_book}, which computes ${{\left( \mathbf{A}_{p}^{m+1} \right)}^{+ }}={{\left( \mathbf{A}_{{}}^{m+1} \right)}^{+ }}$ from ${{\left( \mathbf{A}_{0}^{m+1} \right)}^{+ }}={{\left( \mathbf{A}_{{}}^{m} \right)}^{+ }}$ in $p$ iterations, by (\ref{A_1col_inv_book1}),
%(\ref{d_col_define4311}),
%(\ref{c_col_define1})
%and
%(\ref{b_col_def1ab}).

\begin{algorithm}%[H]
\caption{Pseudoinverse of Col.-Partitioned
Matrix}%Ëã·¨Ãû×Ö
\LinesNumbered %ÒªÇóÏÔʾÐкÅ
\KwIn{$\mathbf{H}_{m+1}$ with $p$ columns, $\mathbf{A}_{0}^{m+1}=\mathbf{A}^{m}$, ${{(\mathbf{A}_{0}^{m+1})}^{+}}$, the ridge parameter $\lambda \to 0$},   % ÊäÈë²ÎÊý
\KwOut{${{({{\mathbf{A}}_p^{m+1}})}^{+}}={{({{\mathbf{A}}^{m+1}})}^{+}}={{\left[ {{\mathbf{A}}^{m}}|{{\mathbf{H}}_{m+1}} \right]}^{+}}$ }%Êä³ö
$i=0$\; %\;ÓÃÓÚ»»ÐÐ
\While{$i<p$}{$\mathbf{D}={{(\mathbf{A}_{i}^{m+1})}^{+}}\mathbf{H}_{m+1}^{:,i+1:p}$,
$\mathbf{C}= \mathbf{H}_{m+1}^{:,i+1:p} -\mathbf{A}_{i}^{m+1}\mathbf{D}$\;
$[k,{{\mathbf{G}}_{k}}] = {\mathop{\rm InvChol}\nolimits}(\mathbf{C},\lambda)$\;
$\delta =1$\;
\eIf{$k = 0$ (i.e., ${\left|  {{\mathbf{\tilde{c}}}_{:{1}}} \right|^2}< \varepsilon$)}{\For{$j=2:p-i$}{
¡¡¡¡\eIf{$ {\left|  {{\mathbf{{c}}}_{:\text{j}}} \right|^2}< \varepsilon$}{
¡¡¡¡¡¡¡¡$\delta =\delta +1$;
¡¡¡¡}{break;  //Terminate the for loop }}}{$\mathbf{B}_{k}^{T}={{\mathbf{G}}_{k}}\mathbf{G}_{k}^{T}\mathbf{C}_{k}^{T}$\;
${{(\mathbf{A}_{i+k}^{m+1})}^{+}}= \left[ \begin{matrix}
   {{(\mathbf{A}_{i}^{m+1})}^{+}}-{{\mathbf{D}}_{k}}\mathbf{B}_{k}^{T}  \\
   \mathbf{B}_{k}^{T}  \\
\end{matrix} \right]\text{ }$\;
$i=i+k$;}
\If{$i < p$}{$\mathbf{\tilde H}_{\delta}=\mathbf{H}_{m+1}^{:,i+1:i+\delta}$ \;
${{\mathbf{D}}_{\delta }} =
\begin{cases}
\mathbf{D}^{:,1:\delta} \quad if\ k = 0  \\
{{(\mathbf{A}_{i}^{m+1})}^{+}}\mathbf{\tilde H}_{\delta} \quad if\ k \ne 0;
\end{cases}$ \\
${{\bf{B}}_{\delta}^T} =
\begin{cases}
  {{({\bf{I}} + {{\bf{D}}_{\delta }^T}{{\bf{D}}_{\delta }} )^{ - 1}}} {{\bf{\tilde D}}_{\delta }} \quad if \   m \ge max(n+i,{\delta })     \\
{{({\bf{I}} + {{\mathbf{\tilde D}}_{\delta }} \mathbf{\tilde H}_{\delta} )}^{ - 1}} {{\mathbf{\tilde D}}_{\delta }} \quad \quad \quad  if\ n+i \ge m \ge {\delta } \\
{{\mathbf{\tilde D}}_{\delta }}{{({\bf{I}} + \mathbf{\tilde H}_{\delta} {{\mathbf{\tilde D}}_{\delta }} )}^{ - 1}} \quad \quad \quad \quad \quad if\ m \le {\delta },
\end{cases}$
where $\mathbf{\tilde H}_{\delta}$ is $m \times  {\delta}$, and ${{\mathbf{\tilde D}}_{\delta }}$ is computed by
${{\mathbf{\tilde D}}_{\delta }}={{\mathbf{D}}_{\delta}^T}{{(\mathbf{A}_{i}^{m+1})}^{+}}$ if required;    \\
${{(\mathbf{A}_{i+\delta}^{m+1})}^{+}} = \left[ \begin{matrix}
   {{(\mathbf{A}_{i}^{m+1})}^{+}}-{{\mathbf{D}}_{\delta }}\mathbf{B}_{\delta }^{T}  \\
   \mathbf{B}_{\delta }^{T}  \\
\end{matrix} \right]$; \\
$i=i+\delta $;}}
\end{algorithm}

%\begin{algorithm}
%\caption{My ${\mathop{\rm InvChol}\nolimits}$ best}%Ëã·¨Ãû×Ö
%\SetAlgoLined
%\SetKwProg{Fn}{Function}{}{end}
%\Fn{$[k, {{\mathbf{G}}_{k}}]={\mathop{\rm InvChol}\nolimits}(\mathbf{C}, \lambda)$}
%{\For{$k=1:size(\mathbf{C},2)$}{
%¡¡¡¡${{\mathbf{\tilde{c}}}_{:{k}}}={{\mathbf{c}}_{:k}}-{{\mathbf{C}}_{k-1}}{{\mathbf{G}}_{k-1}}\mathbf{G}_{k-1}^{T}\mathbf{C}_{k-1}^{T}{{\mathbf{c}}_{:k}}$ (${{\mathbf{\tilde{c}}}_{:1}}={{\mathbf{c}}_{:1}}$)\;
%\eIf{$ {\left|  {{\mathbf{\tilde{c}}}_{:{k}}} \right|^2}< \varepsilon$}
%{$k=k-1$\;
%break; //Terminate the for loop
%¡¡¡¡}
%{¡¡¡¡¡¡¡¡${{f}_{kk}}=1/\sqrt{\mathbf{\tilde{c}}_{:k}^{T}{{{\mathbf{\tilde{c}}}}_{:k}}}$\;
%${{\mathbf{\tilde{f}}}_{k-1}}=-{{f}_{kk}}{{\mathbf{G}}_{k-1}}\mathbf{G}_{k-1}^{T}\mathbf{C}_{k-1}^{T}{{\mathbf{c}}_{:k}}$\;
%${{\mathbf{G}}_{k}}=\left[ \begin{matrix}
%   {{\mathbf{G}}_{k-1}} & {{{\mathbf{\tilde{f}}}}_{k-1}}  \\
%   \mathbf{0}_{k-1}^{T} & {{f}_{kk}}  \\
%\end{matrix} \right]$ (${{\mathbf{G}}_{1}}=\left[ {{f}_{11}} \right]$);
%}}
%\textbf{return} $k$, ${{\mathbf{G}}_{k}}$\;}
%\end{algorithm}

\section{The Proposed Algorithm for Column-Partitioned Matrices in BLS}

The algorithm for the pseudoinverse of a column-partitioned
matrix is described in \textbf{Algorithm 1}.
 In \textbf{Algorithm 1},  there is a while loop including all $26$ rows except row $1$,
 of which the first iteration will be introduced
%Thus
 in what follows.
% we introduce of this while loop.

\begin{algorithm}
\caption{The ${\mathop{\rm InvChol}\nolimits}$ function implemented with the inverse Cholesky factorization}%Ëã·¨Ãû×Ö
\SetAlgoLined
\SetKwProg{Fn}{Function}{}{end}
\Fn{$[k, {{\mathbf{G}}_{k}}]={\mathop{\rm InvChol}\nolimits}(\mathbf{C}, \lambda)$}
{\For{$k=1:{\mathop{\rm size}\nolimits}(\mathbf{C},2)$}{
¡¡¡¡$\left|  {{\mathbf{\tilde{c}}}_{:{k}}} \right|^2=\mathbf{c}_{:k}^{T}{{\mathbf{c}}_{:k}}-\mathbf{c}_{:k}^{T}{{\mathbf{C}}_{k-1}}{{\mathbf{G}}_{k-1}}\mathbf{G}_{k-1}^{T}\mathbf{C}_{k-1}^{T}{{\mathbf{c}}_{:k}}$ ($\left|  {{\mathbf{\tilde{c}}}_{:{1}}} \right|^2=\mathbf{c}_{:1}^{T}{{\mathbf{c}}_{:1}}$)\;
\eIf{$ {\left|  {{\mathbf{\tilde{c}}}_{:{k}}} \right|^2}< \varepsilon$}
{$k=k-1$\;
break; //Terminate the for loop
¡¡¡¡}
{¡¡¡¡¡¡¡¡${{g}_{kk}}=1/\sqrt{\left|  {{\mathbf{\tilde{c}}}_{:{k}}} \right|^2}$\;
${{\mathbf{\tilde{g}}}_{:k}}=-{{g}_{kk}}{{\mathbf{G}}_{k-1}}\mathbf{G}_{k-1}^{T}\mathbf{C}_{k-1}^{T}{{\mathbf{c}}_{:k}}$\;
${{\mathbf{G}}_{k}}=\left[ \begin{matrix}
   {{\mathbf{G}}_{k-1}} & {{{\mathbf{\tilde{g}}}}_{:k}}  \\
   \mathbf{0}_{k-1}^{T} & {{g}_{kk}}  \\
\end{matrix} \right]$ (${{\mathbf{G}}_{1}}=\left[ {{g}_{11}} \right]$);
}}
\textbf{return} $k$, ${{\mathbf{G}}_{k}}$\;}
\end{algorithm}

%\SetKwFunction{FMain}{\bf{InvChol}}
%\SetKwProg{Fn}{Function}{:}{}
%\Fn{\FMain{}}

%\begin{algorithm}
%\caption{$\bf{InvChol}$ by Chol Matlab}%Ëã·¨Ãû×Ö
%\SetAlgoLined
%\SetKwProg{Fn}{Function}{}{end}
%\Fn{$[k, {{\mathbf{G}}_{k}}]=\mathbf{InvChol}(\mathbf{C}, \lambda)$}
%{$[{{\mathbf{\tilde{G}}}}, Flag]=\mathbf{Chol}(\mathbf{C}^{T}\mathbf{C}+\lambda\mathbf{I})$ \;
%${{\mathbf{\tilde{c}}}_{:{k}}}={{\mathbf{c}}_{:k}}-{{\mathbf{C}}_{k-1}}{{\mathbf{G}}_{k-1}}\mathbf{G}_{k-1}^{T}\mathbf{C}_{k-1}^{T}{{\mathbf{c}}_{:k}}$\;
%\If{$\sim isempty({{\mathbf{\tilde{G}}}})$}{$\mathbf{\tilde g}=diag({\mathbf{\tilde{G}}})$\;
%$ \mathbf{d} = find(\mathbf{\tilde g}< \varepsilon)$\;
%\If{$\sim isempty(\mathbf{d})$}{ $Flag = \mathbf{d}(1)$\;
%${\mathbf{\tilde{G}}}={\mathbf{\tilde{G}}}(1:Flag-1,1:Flag-1)$;}
%\If{$\sim isempty({\mathbf{\tilde{G}}})$}{${{\mathbf{G}}_{k}}={\mathbf{\tilde{G}}^{-1}}$;}}
%$k=size({{\mathbf{G}}_{k}},1)$\;
%\textbf{return} $k$, ${{\mathbf{G}}_{k}}$\;}
%\end{algorithm}

 The index $i$ for the while loop denotes that the pseudoinverse
 $\mathbf{A}_{i}^{m+1}=\left[ {{\mathbf{A}}^{m}}|\mathbf{H}_{m+1}^{:,1:i} \right]$
 (defined by (\ref{A_k_def12433}))
 is available, and  the initial  $i$  is set to $0$ in row 1.
 %Firstly
 In row $3$,
%^and $4$,
 $\mathbf{D}$ and $\mathbf{C}$ are computed by (\ref{D_matrix_def_111}) and (\ref{C_matrix_def_111}), respectively.
 Then in row $4$, the function  $[k, {{\mathbf{G}}_{k}}]={\mathop{\rm InvChol}\nolimits}(\mathbf{C}, \lambda)$
 defined in \textbf{Algorithm 2} is  applied  to   find the minimum $k \ge 0$
%($k=0,1,\cdots$)
 satisfying ${{\mathbf{\tilde{c}}}_{:{k+1}}}={\bf{0}}$ (i.e., $ {\left|  {{\mathbf{\tilde{c}}}_{:{k+1}}} \right|^2}<\varepsilon$ where $\varepsilon \to 0$ is  a positive
number near zero,  e.g.,  $\varepsilon=10^{-10}$)
and the corresponding $k \times k$ inverse Cholesky factor ${{\mathbf{G}}_{k}}$ satisfying (\ref{R_def12445}),
or find the $k$ equal to the column number of $\mathbf{C}$ and the corresponding
% $k \times k$ inverse Cholesky factor
 ${{\mathbf{G}}_{k}}$. In \textbf{Algorithm 2},
 $\left|{{\mathbf{\tilde{c}}}_{:{k}}} \right|^2$ is computed by~\footnote{We can also use (\ref{c_col_define4myChol}) instead of (\ref{c_lentghcCFFcc}), at the cost of higher complexity.} (\ref{c_lentghcCFFcc}) for $k=1,2, \cdots$,
 till the first ${\left|  {{\mathbf{\tilde{c}}}_{:{k}}} \right|^2}=0$ (i.e., ${\left|  {{\mathbf{\tilde{c}}}_{:{k}}} \right|^2}< \varepsilon$)
 or $k$ reaches the column number of $\mathbf{C}$. When ${\left|  {{\mathbf{\tilde{c}}}_{:{k}}} \right|^2} \ne 0$,
  ${{\mathbf{G}}_{k}}$
%required in (\ref{comments_myB_def111})
%can be
is computed iteratively from $\mathbf{C}_{k}$
 by
%(\ref{c_lentghcCFFcc}),
(\ref{ita_myChol4c_wave}), (\ref{u_myChol_12321}) and (\ref{F_iterative_Add}).
%  which are summarized in the function  $[k, {{\mathbf{G}}_{k}}]={\mathop{\rm InvChol}\nolimits}(\mathbf{C}, \lambda)$
% defined in \textbf{Algorithm 2}.
Notice that the above function
 $[k, {{\mathbf{G}}_{k}}]={\mathop{\rm InvChol}\nolimits}(\mathbf{C}, \lambda)$
 can also be implemented with  \textbf{Algorithm 3} instead of
 \textbf{Algorithm 2}, when the Matlab built-in function ``chol" is preferred.
% , we can use, to define
% utilized in row  $4$ of \textbf{Algorithm 1}.
In \textbf{Algorithm 2} and \textbf{Algorithm 3},
%Moreover,
the positive real number
$\lambda$ is the ridge parameter satisfying  $\lambda \to 0$, which
%is utilized in
% the generalized inverse with the ridge regression approximation, where
%  the ridge parameter $\lambda \to 0$
   is utilized to approximate
   the generalized inverse with the ridge inverse~\cite{best_ridge_inv_paper213}, as in the original BLS~\cite{BL_trans_paper}.

 If $k = 0$, i.e., $\mathbf{\tilde{c}}_{:k+1}^{{}}=\mathbf{\tilde{c}}_{:1}^{{}}= \mathbf{0}$, ${\delta }$ is decided in rows $7-13$, which means that the first ${\delta }$ columns of $\mathbf{C}$ are zeros; Otherwise in rows $15-17$, the $k \times k$ matrix ${{\mathbf{G}}_{k}}$  is applied to compute $\mathbf{B}_{k}^{T}$ by (\ref{comments_myB_def111}), and then $\mathbf{B}_{k}^{T}$  is applied to
 to update ${{(\mathbf{A}_{i}^{m+1})}^{+}}$ into ${{(\mathbf{A}_{i+k}^{m+1})}^{+}}$ by (\ref{BLK_general_inv_def1223}).
 Moreover, if $i<p$ after the above operations,   ${{\mathbf{\tilde{c}}}_{:{k+1}}}={\bf{0}}$ must have been found
 % find the minimum $k \ge 0$
%%($k=0,1,\cdots$)
% satisfying
  in row $4$, and the first ${\delta } \ge 1$ columns of $\mathbf{C}$ are zeros when $k=0$.
   Thus the the  pseudoinverse
 ${{(\mathbf{A}_{i}^{m+1})}^{+}}$ is updated into ${{(\mathbf{A}_{i+\delta}^{m+1})}^{+}}$   by (\ref{BLK_general_inv_def1223})
   in rows $20-24$, where $\mathbf{B}_{\delta }^{T}$ is computed by  (\ref{B_Matrix_def2bNEWColRankOriginalPosDef})/(\ref{B_Matrix_def2bNEWColRankOtherPosDef})/(\ref{B_Matrix_def2bNEWColRankPosDef}),
   %(\ref{B_Matrix_def2bOld})/(\ref{B_Matrix_def2bOther})/(\ref{B_Matrix_def2b}),???
and   ${{\mathbf{D}}_{\delta }}=\mathbf{D}^{:,1:\delta}$ is the first $\delta$ columns of ${{\mathbf{D}}}$ if $k = 0$,  or %${\delta }=1$, while
 %${{\mathbf{D}}_{\delta }}={{\mathbf{D}}_{1}}$ is
  computed by (\ref{DjMinus1EquAH})  % D_matrix_def_111
  if $k \ne 0$.

\begin{algorithm}
\caption{The ${\mathop{\rm InvChol}\nolimits}$ function implemented with the Matlab built-in function ``chol"}%Ëã·¨Ãû×Ö
\SetAlgoLined
\SetKwProg{Fn}{Function}{}{end}
\Fn{$[k, {{\mathbf{G}}_{k}}]={\mathop{\rm InvChol}\nolimits}(\mathbf{C}, \lambda)$}
{$[{{\mathbf{\tilde{G}}}},FLAG]={\mathop{\rm chol}\nolimits}(\mathbf{C}^{T}\mathbf{C}+\lambda\mathbf{I})$ \;
\If{$\sim {\mathop{\rm isempty}\nolimits}({{\mathbf{\tilde{G}}}})$}{$\mathbf{\tilde g}={\mathop{\rm diag}\nolimits}({\mathbf{\tilde{G}}})$,
$ \mathbf{d} = {\mathop{\rm find}\nolimits}(\mathbf{\tilde g}< \varepsilon)$\;
\If{$\sim {\mathop{\rm isempty}\nolimits}(\mathbf{d})$}{${\mathbf{\tilde{G}}}={\mathbf{\tilde{G}}}(1:\mathbf{d}(1)-1,1:\mathbf{d}(1)-1)$;}}
${{\mathbf{G}}_{k}}={\mathbf{\tilde{G}}^{-1}}$, $k= {\mathop{\rm size}\nolimits}({{\mathbf{G}}_{k}},1)$\;
\textbf{return} $k$, ${{\mathbf{G}}_{k}}$\;}
\end{algorithm}
%$Flag = \mathbf{d}(1)$\;

  %In fact
  In the  first iteration,    (\ref{B_Matrix_def2aNEWColRankPosDef}) is implemented in row $15$ if
 $k=p$,
%$k=p \ge 1$ in row $17$,
%and
% actually
 %(\ref{B_Matrix_def2b})/(\ref{B_Matrix_def2bOther})
% (\ref{B_Matrix_def2bOld}), (\ref{B_Matrix_def2bOther}) and (\ref{B_Matrix_def2b})
%(\ref{B_Matrix_def2bNEWColRankOriginalPosDef}), (\ref{B_Matrix_def2bNEWColRankOtherPosDef})
%and (\ref{B_Matrix_def2bNEWColRankPosDef})
 (\ref{B_Matrix_def2bNEWColRankOriginalPosDef})/(\ref{B_Matrix_def2bNEWColRankOtherPosDef})/(\ref{B_Matrix_def2bNEWColRankPosDef})
  is implemented in row $22$ if $\delta=p$, and (\ref{B_Matrix_def2cNEWColRankPosDef})
 %Moreover, actually
 is corresponding to
 all other cases. Moreover, if $i<p$ after the above-described first iteration of the while loop,
 %and correspondingly
 the next iteration of the while loop
 %is started
 will start with $\mathbf{H}_{m+1}^{:,i+1:p}$ (including only the last  $p-i$ columns
 of ${{\mathbf{H}}_{m+1}}$),  ${{\mathbf{A}_{i}^{m+1}}}=\left[ {{\mathbf{A}}^{m}}|\mathbf{H}_{m+1}^{:,1:i} \right]$  and
 ${{(\mathbf{A}_{i}^{m+1})}^{+}}$.

\section{The Proposed Algorithm for Row-Partitioned Matrices in BLS}
%${{{{\mathbf{A}}^{+ }}}}=\left({{({{\mathbf{A}}^{T}})}^{+}}\right)^{T}$

The incremental learning for the increment of input data
in \cite{BL_trans_paper} utilizes
the pseudoinverse of
 the  row-partitioned matrix
 %that should read  \cite{my_correction2BL}
\begin{equation}\label{AxInputIncrease31232USE79547}
{}^x{\bf{A}}_n^m = \left[ \begin{array}{l}
{\bf{A}}_n^m\\
{\bf{A}}_x^{}
\end{array} \right],
\end{equation}%[AxInputIncrease31232]
where ${\bf{A}}_n^m$ is
$m \times n$, and ${\bf{A}}_x^{}$ can be assumed to be
$q \times n$.
%$(m+q) \times n$ ${}^x{\bf{A}}_n^m$
%We need to
%Let us use
Equation (c) in \cite[Ex. 1.16]{cite_general_inv_book} can be written as
 \begin{equation}\label{AinvT2AtInv34985}
{{({{\mathbf{A}}^{+ }})}^{T}}={{({{\mathbf{A}}^{T}})}^{+}},
\end{equation}
into which we can substitute (\ref{AxInputIncrease31232USE79547})
%(\ref{AinvT2AtInv34985})
to obtain
 \begin{equation}\label{AinvT2AtInv34985Other1}
{{\left({{({}^x{\bf{A}}_n^m)}^{+ }}\right)}^{T}}={{\left[ {{{({\bf{A}}_n^m)}^T}|{\bf{A}}_x^T} \right]}^{+}}.
\end{equation}
Then substitute (\ref{AinvT2AtInv34985Other1})
into
(\ref{A_1col_inv_book1matrix})
to obtain
 \begin{equation}\label{A_1col_inv_book1matrixOther1}
 {{\left({{({}^x{\bf{A}}_n^m)}^{+ }}\right)}^{T}}
=\left[ \begin{matrix}
   {{\left({{({\bf{A}}_n^m)}^T}\right)}^{+}}-\mathbf{D}{{\mathbf{B}}^{T}}  \\
   {{\mathbf{B}}^{T}}  \\
\end{matrix} \right],
\end{equation}
%which can be written as,
i.e.,
 \begin{equation}\label{A_1col_inv_book1matrixOther222}
 {{{({}^x{\bf{A}}_n^m)}^{+ }}}
=\left[ {{({\bf{A}}_n^m)}^{+}}-{\mathbf{B}} {{\mathbf{D}}^{T}}  | {\mathbf{B}}  \right].
\end{equation}

%Pseudoinverse of Row-Partitioned Matrix
\begin{algorithm}%[H]
\caption{Pseudoinverse of Row-Partitioned
Matrix}%Ëã·¨Ãû×Ö
\LinesNumbered %ÒªÇóÏÔʾÐкÅ
\KwIn{${\bf{A}}_x^{}$ with $q$ rows, ${}_0^x{\bf{A}}_n^m={\bf{A}}_n^m$, ${{({}_0^x{\bf{A}}_n^m)}^{+}}$, the ridge parameter $\lambda \to 0$}% ÊäÈë²ÎÊý
\KwOut{${{({}_q^x{\bf{A}}_n^m)}^{+}}={{({}^x{\bf{A}}_n^m)}^{+}}={{\left[ \begin{array}{l}
{\bf{A}}_n^m\\
{\bf{A}}_x
\end{array} \right]}^{+}}$ }%Êä³ö
$i=0$\; %\;ÓÃÓÚ»»ÐÐ
\While{$i<q$}{${{\mathbf{D}}^{T}} = {\bf{A}}_x^{i+1:q,:} {{({}_i^x{\bf{A}}_n^m)}^{+}}$\;
$\mathbf{C}={ ({\bf{A}}_x^{i+1:q,:})^T}-{{({}_i^x{\bf{A}}_n^m)}^T}\mathbf{D}$\;
$[k,{{\mathbf{G}}_{k}}] = {\mathop{\rm InvChol}\nolimits}(\mathbf{C},\lambda)$\;
$\delta =1$\;
\eIf{$k = 0$ (i.e., ${\left|  {{\mathbf{\tilde{c}}}_{:{1}}} \right|^2}< \varepsilon$)}{\For{$j=2:q-i$}{
¡¡¡¡\eIf{$ {\left|  {{\mathbf{{c}}}_{:\text{j}}} \right|^2}<\varepsilon$}{
¡¡¡¡¡¡¡¡$\delta =\delta +1$;
¡¡¡¡}{break;  //Terminate the for loop }}}{${{\bf{B}}_{k}} ={{\mathbf{C}}_{k}} {{\mathbf{G}}_{k}}\mathbf{G}_{k}^{T}$\;
${\left({{({}_{i}^x{\bf{A}}_n^m)}^{+ }}\right)}
=\left[ {{({}_{i-k}^{\ \  x}{\bf{A}}_n^m)}^{+}}-{\mathbf{B}}_{k} {{\mathbf{D}}_{k}^{T}}  | {\mathbf{B}}_{k}  \right]$;
$i=i+k$;}
\If{$k \le p-i-1$}{${\bf{A}}_x^{\delta}={\bf{A}}_x^{i+1:i+\delta,:}$\;
${{\mathbf{D}}_{\delta}^{T}} =
\begin{cases}
(\mathbf{D}^{:,1:\delta})^{T} \quad if\ k = 0  \\
 {\bf{A}}_x^{\delta}{{({}_{i-\delta}^{\ \  x}{\bf{A}}_n^m)}^{+}} \quad if\ k \ne 0
\end{cases}$\;
${{\bf{B}}_{\delta}^T} =
\begin{cases}
{\bf{\tilde D}}^T {{({\bf{I}} + {{\bf{D}}^T}{\bf{D}} )^{ - 1}}}\ if \,  n \ge max(m+i,\delta)   \\
{{\mathbf{D}}_{\delta}^{T}}{{({\bf{I}} + {\bf{A}}_x^{\delta}{{\mathbf{D}}_{\delta}^{T}})}^{ - 1}} \quad if\ m+i \ge n \ge \delta \\
{{({\bf{I}} + {{\mathbf{D}}_{\delta}^{T}} {\bf{A}}_x^{\delta} )}^{ - 1}}{{\mathbf{D}}_{\delta}^{T}} \quad if\  n  \le {\delta }
\end{cases}$
%{\bf{\tilde D}}^T {{({\bf{I}} + {{\bf{D}}^T}{\bf{D}} )^{ - 1}}}\ if \,  {\bf{C}} = {\bf{0}}, n \ge max(m,q)   &   \label{B_Matrix_def2bFullColumn222aaaorig}  \\
%{\bf{\tilde D}}^T{{({\bf{I}} + {\bf{A}}_x{\bf{\tilde D}}^T)}^{ - 1}} \quad \ if\ {\bf{C}} = {\bf{0}},  m \ge n \ge q &   \label{B_Matrix_def2bFullColumn222bbb}\\
%{{({\bf{I}} + {\bf{\tilde D}}^T {\bf{A}}_x)}^{ - 1}}{\bf{\tilde D}}^T \quad \quad \quad \ if\ {\bf{C}} = {\bf{0}}, n \le q    &  \label{B_Matrix_def2bFullColumn222aaa}\\
where ${\bf{A}}_x^{\delta}$ is ${\delta } \times  n$\;
$ {{{({}_{i+\delta}^{\ \ x}{\bf{A}}_n^m)}^{+ }}}
=\left[ {{({}_{i}^{x}{\bf{A}}_n^m)}^{+}}-{\mathbf{B}}_{\delta } {{\mathbf{D}}_{\delta }^{T}}  | {\mathbf{B}}_{\delta}  \right]$\;
$i=i+\delta$;}}
\end{algorithm}

%From
%(\ref{A_1col_inv_book1matrix})
%%(\ref{AinvT2AtInv34985Other1})
% and (\ref{AinvT2AtInv34985Other1}),  it can be seen that
Obviously
 ${{\left[ {{\mathbf{A}}^{m}}|{{\mathbf{H}}_{m+1}} \right]}}$
in
(\ref{A_1col_inv_book1matrix}) is replaced with ${{\left[ {{{({\bf{A}}_n^m)}^T}|{\bf{A}}_x^T} \right]}}$ in
%(\ref{AinvT2AtInv34985Other1})
(\ref{AinvT2AtInv34985Other1}).
%and then
Accordingly
 in %()()()
  (\ref{D_matrix_def_111}),
(\ref{C_matrix_def_111}), (\ref{DTilde2DtAmInv94832})
and
(\ref{B_Matrix_def2abc}),
  ${{\mathbf{A}}^{m}}$ and ${{\mathbf{H}}_{m+1}}$
  should be replaced by ${{({\bf{A}}_n^m)}^T}$  and ${\bf{A}}_x^T$, respectively, to obtain
% (no index)
 % ,
%where
%\begin{equation}\label{D_matrix_def_111Other1}
$\mathbf{D} = {{\left({{({\bf{A}}_n^m)}^T}\right)}^{+}} {{\bf{A}}_x^T}$
%\end{equation}
%which can be written as
that can be written as
%is
\begin{equation}\label{D_matrix_def_111Other222}
{{\mathbf{D}}^{T}} = {{\bf{A}}_x} {{({\bf{A}}_n^m)}^{+}},
\end{equation}
%\begin{subnumcases}{\label{B_Matrix_def2abcOther1}}
%{{{\bf{B}}^T} = {{({{\mathbf{C}}^{T}}\mathbf{C})}^{-1}}{{\mathbf{C}}^{T}}\quad if\ each\ {\bf{\tilde c}}_{:k}^{} \ne {\bf{0}}(1 \le k \le p)} &  \label{B_Matrix_def2aOther1}\\
%{{{\bf{B}}^T} = {{({\bf{I}} + {{\bf{D}}^T}{\bf{D}})}^{ - 1}}{{\bf{D}}^T} {{\left({{({\bf{A}}_n^m)}^T}\right)}^{+}}  \quad if\ {\bf{C}} = {\bf{0}}}  &  \label{B_Matrix_def2bOther1}\\
%{{{\bf{B}}^T} = \cdots \quad if\ {\bf{C}} \ne {\bf{0}} \ but\ several\ {\bf{\tilde c}}_{:k}^{} = {\bf{0}}}, & \label{B_Matrix_def2cOther1}
%\end{subnumcases}
%\begin{equation}\label{C_matrix_def_111Other1}
%\mathbf{C}={{\bf{A}}_x^T}-{{({\bf{A}}_n^m)}^T}\mathbf{D},
%\end{equation}
\begin{equation}\label{C_matrix_def_111Other222}
\mathbf{C}={{\bf{A}}_x^T}-{{({\bf{A}}_n^m)}^T}\mathbf{D},
\end{equation}
${\bf{\tilde D}}={{\bf{D}}^T}{{({{({\bf{A}}_n^m)}^+})}^T}$ that can be written as
\begin{equation}\label{DTilde2DtAmInv94832Input}
{\bf{\tilde D}}^T={{{{({\bf{A}}_n^m)}^+}}}{{\bf{D}}},
\end{equation}
and
%20190704  20190704 20190704 20190704 20190704 20190704 20190704
%\begin{subequations}{\label{B_Matrix_def2abcOther1}}
% \begin{numcases}
%{  {{\bf{B}}^T} = }
%{{{({{\mathbf{C}}^{T}}\mathbf{C})}^{-1}}{{\mathbf{C}}^{T}}\quad if\ each\ {\bf{\tilde c}}_{:k}^{} \ne {\bf{0}}(1 \le k \le p)} &  \label{B_Matrix_def2aOther1}\\
%{ {{({\bf{I}} + {{\bf{D}}^T}{\bf{D}})}^{ - 1}}{{\bf{D}}^T} {{\left({{({\bf{A}}_n^m)}^T}\right)}^{+}}  \quad if\ {\bf{C}} = {\bf{0}}}  &  \label{B_Matrix_def2bOther1}\\
%{\cdots \quad if\ {\bf{C}} \ne {\bf{0}} \ but\ several\ {\bf{\tilde c}}_{:k}^{} = {\bf{0}}}, & \label{B_Matrix_def2cOther1}
%\end{numcases}
%\end{subequations}
$$
{{\bf{B}}^T} =
\begin{cases}
  {{({\bf{I}} + {{\bf{D}}^T}{\bf{D}} )^{ - 1}}}
{\bf{\tilde D}} \ if \,  {\bf{C}} = {\bf{0}}, n \ge max(m,q)    \\
{{({\bf{I}} + {\bf{\tilde D}}{\bf{A}}_x^T)}^{ - 1}}{\bf{\tilde D}} \quad \quad if\ {\bf{C}} = {\bf{0}},   m \ge n \ge q  \\
{\bf{\tilde D}}{{({\bf{I}} + {\bf{A}}_x^T{\bf{\tilde D}})}^{ - 1}} \quad \quad \quad \quad if\ {\bf{C}} = {\bf{0}},  n \le q      \\
{{{\bf{C}}^ + }\quad \quad \quad \quad \quad if\ each\ {\bf{\tilde c}}_{:k}^{} \ne {\bf{0}}(1 \le k \le q)} \\
{\cdots \quad \quad \quad \quad  if\ {\bf{C}} \ne {\bf{0}} \ but\ several\ {\bf{\tilde c}}_{:k}^{} = {\bf{0}}}
\end{cases}
$$
%\begin{subequations}{\label{B_Matrix_def2abc}}
% \begin{numcases}
%{  {{\bf{B}}^T} = }
%  {{({\bf{I}} + {{\bf{D}}^T}{\bf{D}} )^{ - 1}}}
%{\bf{\tilde D}} \ if \,  {\bf{C}} = {\bf{0}}, m \ge max(n,p)     &  \label{B_Matrix_def2bOld} \\
%{{({\bf{I}} + {\bf{\tilde D}}{{\mathbf{H}}_{m+1}})}^{ - 1}}{\bf{\tilde D}} \ if\ {\bf{C}} = {\bf{0}}, n \ge m \ge p  &  \label{B_Matrix_def2bOther} \\
%{\bf{\tilde D}}{{({\bf{I}} + {{\mathbf{H}}_{m+1}}{\bf{\tilde D}})}^{ - 1}} \quad \quad \ if\ {\bf{C}} = {\bf{0}},  m \le p  &  \label{B_Matrix_def2b} \\
% {{{\bf{C}}^ + }\quad \quad \quad \quad \ if\ each\ {\bf{\tilde c}}_{:k}^{} \ne {\bf{0}}(1 \le k \le p)}  &  \label{B_Matrix_def2a} \\
%{\cdots \quad \quad \quad \  if\ {\bf{C}} \ne {\bf{0}} \ but\ several\ {\bf{\tilde c}}_{:k}^{} = {\bf{0}}},  &  \label{B_Matrix_def2c}
%\end{numcases}
%\end{subequations}
 % which can be written as
 that  can be written as
 %is
%\begin{subnumcases}{\label{B_Matrix_def2abcOther1}}
%{{{\bf{B}}^T} = {{({{\mathbf{C}}^{T}}\mathbf{C})}^{-1}}{{\mathbf{C}}^{T}}\quad if\ each\ {\bf{\tilde c}}_{:k}^{} \ne {\bf{0}}(k=1,\cdots,q)} &  \label{B_Matrix_def2aOther1}\\
%{{{\bf{B}}^T} = {{({\bf{I}} + {{\bf{D}}^T}{\bf{D}})}^{ - 1}}{{\bf{D}}^T}{{({{\bf{A}}^m})}^ + }\quad if\ {\bf{C}} = {\bf{0}}}  &  \label{B_Matrix_def2bOther1}\\
%{{{\bf{B}}^T} = \cdots \quad if\ {\bf{C}} \ne {\bf{0}} \ but\ several\ {\bf{\tilde c}}_{:k}^{} = {\bf{0}}}, & \label{B_Matrix_def2cOther1}
%\end{subnumcases}
%From ()()() we can deduce
%\begin{subnumcases}{\label{B_Matrix_def2abcOther222}}
%{{{\bf{B}}} ={{\mathbf{C}}} {{({{\mathbf{C}}^{T}}\mathbf{C})}^{-1}}\quad if\ each\ {\bf{\tilde c}}_{:k}^{} \ne {\bf{0}}(1 \le k \le p)} &  \label{B_Matrix_def2aOther222}\\
%{{{\bf{B}}} = {{({\bf{A}}_n^m)}^{+}} {{\bf{D}}} {{({\bf{I}} + {{\bf{D}}^T}{\bf{D}})}^{ - 1}}  \quad if\ {\bf{C}} = {\bf{0}}}  &  \label{B_Matrix_def2bOther222}\\
%{{{\bf{B}}} = \cdots \quad if\ {\bf{C}} \ne {\bf{0}} \ but\ several\ {\bf{\tilde c}}_{:k}^{} = {\bf{0}}}, & \label{B_Matrix_def2cOther222}
%\end{subnumcases}
%and
\begin{subequations}{\label{B_Matrix_def2abcOther222}}
 \begin{numcases}
{  {{\bf{B}}} = }
{\bf{\tilde D}}^T {{({\bf{I}} + {{\bf{D}}^T}{\bf{D}} )^{ - 1}}}\ if \,  {\bf{C}} = {\bf{0}}, n \ge max(m,q)   &   \label{B_Matrix_def2bOther222orig}  \\
{\bf{\tilde D}}^T{{({\bf{I}} + {\bf{A}}_x{\bf{\tilde D}}^T)}^{ - 1}} \quad \ if\ {\bf{C}} = {\bf{0}},  m \ge n \ge q  &   \label{B_Matrix_def2bOther222bbb}\\
{{({\bf{I}} + {\bf{\tilde D}}^T {\bf{A}}_x)}^{ - 1}}{\bf{\tilde D}}^T \quad \quad \quad \ if\ {\bf{C}} = {\bf{0}}, n \le q   &  \label{B_Matrix_def2bOther222aaa}\\
{ ({{\bf{C}}^ + })^T \quad \quad \quad \quad if\ each\ {\bf{\tilde c}}_{:k}^{} \ne {\bf{0}}(1 \le k \le q)} &  \label{B_Matrix_def2aOther222}\\
{ \cdots \quad \quad \quad \quad if\ {\bf{C}} \ne {\bf{0}} \ but\ several\ {\bf{\tilde c}}_{:k}^{} = {\bf{0}}}, & \label{B_Matrix_def2cOther222}
\end{numcases}
\end{subequations}
where ${\bf{A}}_x^{}$ is $q \times n$.
%\begin{subequations}{\label{B_Matrix_def2abcOther222}}
% \begin{numcases}
%{  {{\bf{B}}} = }
%{ ({{\bf{C}}^ + })^T \quad if\ each\ {\bf{\tilde c}}_{:k}^{} \ne {\bf{0}}(1 \le k \le q)} &  \label{B_Matrix_def2aOther222}\\
%{{({\bf{I}} + {\bf{A}}_x^T{\bf{\tilde D}})}^{ - 1}}{\bf{\tilde D}}^T \quad if\ {\bf{C}} = {\bf{0}} \And m \le p  &  \label{B_Matrix_def2bOther222aaa}\\
%{\bf{\tilde D}}^T{{({\bf{I}} + {\bf{\tilde D}}{\bf{A}}_x^T)}^{ - 1}} \quad if\ {\bf{C}} = {\bf{0}}  \And m \ge p &   \label{B_Matrix_def2bOther222bbb}\\
%{ \cdots \quad if\ {\bf{C}} \ne {\bf{0}} \ but\ several\ {\bf{\tilde c}}_{:k}^{} = {\bf{0}}}. & \label{B_Matrix_def2cOther222}
%\end{numcases}
%\end{subequations}
%As ${\bf{\tilde c}}_{:k}^{}$ defined by , ${\bf{\tilde c}}_{:k}^{}$ in (\ref{B_Matrix_def2abcOther222}) can be defined by
In (\ref{B_Matrix_def2abcOther222}), ${\bf{\tilde c}}_{:1}^{}$  can be obtained by (\ref{c1_wave_c1_equal}), and according to Theorem 5,
when ${{\mathbf{\tilde{c}}}_{:i}}\ne 0$ for all $1\le i\le k-1$, ${{\mathbf{\tilde{c}}}_{:k}}$ ($k=2,3,\cdots, p$) can be computed by
(\ref{c_col_define4myChol}), where ${{\mathbf{G}}_{k-1}}$ is the inverse Cholesky factor of $\mathbf{C}_{k-1}^{T}\mathbf{C}_{k-1}$.
${{\mathbf{G}}_{k-1}}$  can be computed by (\ref{ita_myChol_1212}),
(\ref{u_myChol_12321}) and (\ref{F_iterative_Add}) when $k\ge3$, or
 %computed
 by (\ref{initialF1compute}) when $k=2$.

Moreover,  since (\ref{B_Matrix_def2abc}) can be written as (\ref{Bt3Choice98476ColRankPosDef}),  (\ref{B_Matrix_def2abcOther222}) can  be written as
%\begin{subequations}{\label{B_Matrix_def2abcFullColumn974}}
% \begin{numcases}
%{  {{\bf{B}}} = }
%{\bf{\tilde D}}^T {{({\bf{I}} + {{\bf{D}}^T}{\bf{D}} )^{ - 1}}}\ if \,  {\bf{C}} = {\bf{0}}, n \ge max(m,q)   &   \label{B_Matrix_def2bFullColumn222aaaorig}  \\
%{\bf{\tilde D}}^T{{({\bf{I}} + {\bf{A}}_x{\bf{\tilde D}}^T)}^{ - 1}} \quad \ if\ {\bf{C}} = {\bf{0}},  m \ge n \ge q &   \label{B_Matrix_def2bFullColumn222bbb}\\
%{{({\bf{I}} + {\bf{\tilde D}}^T {\bf{A}}_x)}^{ - 1}}{\bf{\tilde D}}^T \quad \quad \quad \ if\ {\bf{C}} = {\bf{0}}, n \le q    &  \label{B_Matrix_def2bFullColumn222aaa}\\
%{  ({{\bf{C}}^ + })^T  \quad \quad \quad  if\ \mathbf{C} \ is \  full \ column \ rank} &  \label{B_Matrix_def2aFullColumn222}\\
%{ \cdots \quad if\ {\bf{C}} \ne {\bf{0}} \ is\ not \  full \ column \ rank}. & \label{B_Matrix_def2cFullColumn222}
%\end{numcases}
%\end{subequations}
\begin{subequations}{\label{B_Matrix_def2abcFullColumn974}}
 \begin{numcases}
{  {{\bf{B}}} = }
{\bf{\tilde D}}^T {{({\bf{I}} + {{\bf{D}}^T}{\bf{D}} )^{ - 1}}}\ if \,  {\bf{C}} = {\bf{0}}, n \ge max(m,q)   &   \label{B_Matrix_def2bFullColumn222aaaorig}  \\
{\bf{\tilde D}}^T{{({\bf{I}} + {\bf{A}}_x{\bf{\tilde D}}^T)}^{ - 1}} \quad \ if\ {\bf{C}} = {\bf{0}},  m \ge n \ge q &   \label{B_Matrix_def2bFullColumn222bbb}\\
{{({\bf{I}} + {\bf{\tilde D}}^T {\bf{A}}_x)}^{ - 1}}{\bf{\tilde D}}^T \quad \quad \quad \ if\ {\bf{C}} = {\bf{0}}, n \le q    &  \label{B_Matrix_def2bFullColumn222aaa}\\
{  ({{\bf{C}}^ + })^T  \quad  \quad   \  if\ \mathbf{C}_{{}}^{T}\mathbf{C} \ is \  positive\  definite} &  \label{B_Matrix_def2aFullColumn222}\\
\cdots \cdots  \quad \quad \quad \begin{array}{l}
 \ if\ {{\bf{C}} \ne 0}\ and \ \mathbf{C}_{{}}^{T}\mathbf{C} \ is  \\
\ not \  positive\  definite
\end{array}. & \label{B_Matrix_def2cFullColumn222}
\end{numcases}
\end{subequations}

%{\cdots \quad  \quad \ if\ \mathbf{C}_{{}}^{T}\mathbf{C} \ is\ not \  positive\  definite}

Obviously, (\ref{B_Matrix_def2abcOther222}) can also  be written as the form that is similar to
(\ref{Bt3Choice98476ColRank}), which is omitted for simplicity.
%${{\bf{B}}} =$
%\begin{subnumcases}{\label{B_Matrix_def2abcOther222}}
%{{{\bf{B}}} ={{\mathbf{C}}} {{({{\mathbf{C}}^{T}}\mathbf{C})}^{-1}}\quad if\ each\ {\bf{\tilde c}}_{:k}^{} \ne {\bf{0}}(1 \le k \le p)} &  \label{B_Matrix_def2aOther222}\\
%{{{\bf{B}}} = {{({\bf{A}}_n^m)}^{+}} {{\bf{D}}} {{({\bf{I}} + {{\bf{D}}^T}{\bf{D}})}^{ - 1}}  \quad if\ {\bf{C}} = {\bf{0}}}  &  \label{B_Matrix_def2bOther222}\\
%{{{\bf{B}}} = \cdots \quad if\ {\bf{C}} \ne {\bf{0}} \ but\ several\ {\bf{\tilde c}}_{:k}^{} = {\bf{0}}}, & \label{B_Matrix_def2cOther222}
%\end{subnumcases}
%\begin{equation}\label{C_matrix_def_111Other222}
%\mathbf{C}={{\bf{A}}_x^T}-{{({\bf{A}}_n^m)}^T}\mathbf{D},
%\end{equation}
%\begin{subnumcases}{\label{B_Matrix_def2abcOther222}}
%{{\bf{B}}}=
%\begin{cases}
%({{\bf{C}}^+})^T={{\mathbf{C}}} {{({{\mathbf{C}}^{T}}\mathbf{C})}^{-1}}  \ \  \quad if \ {{\bf{C}} \ne 0}  &  \label{B_Matrix_def2aOther222} \\
%{  {{({\bf{A}}_n^m)}^{+}} {{\bf{D}}} {{({\bf{I}} + {{\bf{D}}^T}{\bf{D}})}^{ - 1}} \quad if \ {\bf{C}} = 0}  &  \label{B_Matrix_def2bOther222}
%\end{cases}
%\end{subnumcases}

%\section{Input Add the whole algorithm}
Let
\begin{equation}\label{AxInput_k_def76036}
{}_k^x{\bf{A}}_n^m = \left[ \begin{array}{l}
{\bf{A}}_n^m\\
{\bf{A}}_x^{1:k,:}
\end{array} \right],
\end{equation}
 where ${\bf{A}}_x^{1:k£¬:}$ denotes the first $k$ rows of ${\bf{A}}_x^{}$. When $k=0$, ${\bf{A}}_x^{1:k£¬:}$ becomes empty and then (\ref{AxInput_k_def76036}) becomes
\begin{equation}\label{x0AisAmn}
{}_0^x{\bf{A}}_n^m={\bf{A}}_n^m.
\end{equation}
Then the algorithm for the pseudoinverse of a row-partitioned
matrix is shown in \textbf{Algorithm 4},
where $\mathbf{C}_{k}$ and
${\mathbf{D}}_{k}$ denote the first $k$ columns of $\mathbf{C}$ and
${\mathbf{D}}$, respectively, and the function
 $[k, {{\mathbf{G}}_{k}}]={\mathop{\rm InvChol}\nolimits}(\mathbf{C}, \lambda)$
 in row $5$
 can be implemented with  \textbf{Algorithm 2} or  \textbf{Algorithm 3}.

\section{Conclusions}

In BLS,
Greville's method \cite{cite_general_inv_book} has been utilized to propose an effective and efficient
incremental learning system without
retraining the whole network from the beginning.
For a column-partitioned matrix
%\begin{equation}\label{Abig_original}
$\mathbf{A}_{{}}^{m+1}=\left[ {{\mathbf{A}}^{m}}|\mathbf{H}_{m+1}^{{}} \right]$
%\end{equation}
where
%the $m \times p$
 $\mathbf{H}_{m+1}^{{}}$ includes $p$ columns,
Greville's method spends  $p$ iterations to
compute ${{\left( \mathbf{A}_{{}}^{m+1} \right)}^{+ }}$ from ${{\left( \mathbf{A}_{{}}^{m} \right)}^{+ }}$,
where ${{\mathbf{A}}^{+ }}$ denotes the pseudoinverse of the matrix $\mathbf{A}$.  However, the incremental algorithms in \cite{BL_trans_paper}
extend Greville's method
   to compute ${{\left( \mathbf{A}_{{}}^{m+1} \right)}^{+ }}$ from ${{\left( \mathbf{A}_{{}}^{m} \right)}^{+ }}$ by just $1$ iteration,
   which have neglected some possible cases, and need further improvements in efficiency and numerical stability.
   In this paper,
 we propose an efficient and numerical stable algorithm from Greville's method,
 to compute ${{\left( \mathbf{A}_{{}}^{m+1} \right)}^{+ }}$ from ${{\left( \mathbf{A}_{{}}^{m} \right)}^{+ }}$ by just $1$ iteration,
where all possible cases are considered, and  the efficient inverse Cholesky factorization in \cite{my_inv_chol_paper}
 can be applied to further reduce the computational complexity. Finally, we give the whole algorithm for column-partitioned matrices in BLS.
 On the other hand, we also give the proposed algorithm for row-partitioned matrices in BLS.

\appendices
\section{Proof of Theorem 5}
Firstly, we verify (\ref{c_col_define4myChol}).
From (\ref{C_matrix_def_onlyH121}) we deduce that the $k$-th column of $\mathbf{C}$ is
${{\mathbf{c}}_{:k}}={{\mathbf{h}}_{:k}}-{{\mathbf{A}}^{m}}{{({{\mathbf{A}}^{m}})}^{+}}{{\mathbf{h}}_{:k}}$, which is substituted into (\ref{c_col_define4myChol}) to obtain
${{\mathbf{\tilde{c}}}_{:k}}={{\mathbf{h}}_{:k}}-{{\mathbf{A}}^{m}}{{({{\mathbf{A}}^{m}})}^{+}}
{{\mathbf{h}}_{:k}}-{{\mathbf{C}}_{k-1}}{{\mathbf{G}}_{k-1}}\mathbf{G}_{k-1}^{T}\mathbf{C}_{k-1}^{T}({{\mathbf{h}}_{:k}}-{{\mathbf{A}}^{m}}{{({{\mathbf{A}}^{m}})}^{+}}{{\mathbf{h}}_{:k}})$ , i.e.,
\begin{multline}\label{comments_c_wave111}
{{\mathbf{\tilde{c}}}_{:k}}=( {{\mathbf{I}}_{}}-{{\mathbf{A}}^{m}}{{({{\mathbf{A}}^{m}})}^{+}}-{{\mathbf{C}}_{k-1}}{{\mathbf{G}}_{k-1}}\mathbf{G}_{k-1}^{T}\mathbf{C}_{k-1}^{T}
\\
+{{\mathbf{C}}_{k-1}}{{\mathbf{G}}_{k-1}}\mathbf{G}_{k-1}^{T}\mathbf{C}_{k-1}^{T}{{\mathbf{A}}^{m}}{{({{\mathbf{A}}^{m}})}^{+}} ){{\mathbf{h}}_{:k}}.
\end{multline}
Equation (\ref{comments_c_wave111}) can be written as
\begin{multline}\label{c_wave_defIAACFFCh}{{\mathbf{\tilde{c}}}_{:k}}= \\
\left( {{\mathbf{I}}_{}}-{{\mathbf{A}}^{m}}{{({{\mathbf{A}}^{m}})}^{+}}-{{\mathbf{C}}_{k-1}}{{\mathbf{G}}_{k-1}}\mathbf{G}_{k-1}^{T}\mathbf{C}_{k-1}^{T} \right){{\mathbf{h}}_{:k}},
\end{multline}
since $\mathbf{C}_{k-1}^{T}{{\mathbf{A}}^{m}}{{({{\mathbf{A}}^{m}})}^{+}}$ in the 2nd row of (\ref{comments_c_wave111}) always satisfies
%In (\ref{comments_c_wave111}),
\begin{equation}\label{CAA_zeros122345}
\mathbf{C}_{k-1}^{T}{{\mathbf{A}}^{m}}{{({{\mathbf{A}}^{m}})}^{+}}=\mathbf{0},
\end{equation}
% is zero, which will be verified as follows.
which will be verified in the next paragraph.

To verify (\ref{CAA_zeros122345}), write
%From
 (\ref{C_matrix_def_onlyH121}) as
   \begin{displaymath}
{{\mathbf{C}}_{k-1}}=\mathbf{H}_{m+1}^{1:k-1}-\mathbf{A}_{{}}^{m}{{(\mathbf{A}_{{}}^{m})}^{+}}\mathbf{H}_{m+1}^{1:k-1},
 \end{displaymath}
%we deduce that the first $k-1$ columns of $\mathbf{C}$ are
%${{\mathbf{C}}_{k-1}}=\mathbf{H}_{m+1}^{1:k-1}-\mathbf{A}_{{}}^{m}{{(\mathbf{A}_{{}}^{m})}^{+}}\mathbf{H}_{m+1}^{1:k-1}$,
 which is substituted into (\ref{CAA_zeros122345}) to obtain
%${{\left( \mathbf{H}_{m+1}^{1:k-1}-\mathbf{A}_{{}}^{m}{{(\mathbf{A}_{{}}^{m})}^{+}}\mathbf{H}_{m+1}^{1:k-1} \right)}^{T}}{{\mathbf{A}}^{m}}{{({{\mathbf{A}}^{m}})}^{+}}=\mathbf{0}$,  i.e.,
 \begin{multline}\label{HAAminusHAAAA13411}
{{(\mathbf{H}_{m+1}^{1:k-1})}^{T}}{{\mathbf{A}}^{m}}{{({{\mathbf{A}}^{m}})}^{+}}- \\
{{(\mathbf{H}_{m+1}^{1:k-1})}^{T}}{{\left( \mathbf{A}_{{}}^{m}{{(\mathbf{A}_{{}}^{m})}^{+}} \right)}^{T}}{{\mathbf{A}}^{m}}{{({{\mathbf{A}}^{m}})}^{+}}=\mathbf{0}.
 \end{multline}
%${{(\mathbf{H}_{m+1}^{1:k-1})}^{T}}{{\mathbf{A}}^{m}}{{({{\mathbf{A}}^{m}})}^{+}}-{{(\mathbf{H}_{m+1}^{1:k-1})}^{T}}{{\left( \mathbf{A}_{{}}^{m}{{(\mathbf{A}_{{}}^{m})}^{+}} \right)}^{T}}{{\mathbf{A}}^{m}}{{({{\mathbf{A}}^{m}})}^{+}}=\mathbf{0}$,
 %into which
 Substitute  (\ref{MoorePenroseDef121c}) into (\ref{HAAminusHAAAA13411}) to obtain
 %write (\ref{HAAminusHAAAA13411}) as
   \begin{multline}
{{(\mathbf{H}_{m+1}^{:,i+1:p})}^{T}}{{\mathbf{A}}^{m}}{{({{\mathbf{A}}^{m}})}^{+}}- \\
 {{(\mathbf{H}_{m+1}^{:,i+1:p})}^{T}}
{{\mathbf{A}}^{m}}{{({{\mathbf{A}}^{m}})}^{+}}{{\mathbf{A}}^{m}}{{({{\mathbf{A}}^{m}})}^{+}}=\mathbf{0},
 \end{multline}
%${{(\mathbf{H}_{m+1}^{:,i+1:p})}^{T}}{{\mathbf{A}}^{m}}{{({{\mathbf{A}}^{m}})}^{+}}-{{(\mathbf{H}_{m+1}^{:,i+1:p})}^{T}}
%{{\mathbf{A}}^{m}}{{({{\mathbf{A}}^{m}})}^{+}}{{\mathbf{A}}^{m}}{{({{\mathbf{A}}^{m}})}^{+}}=\mathbf{0}$,
into which substitute (\ref{MoorePenroseDef121a}) to obtain
%write it as
\begin{equation}\label{HAAminusHAA}
{{(\mathbf{H}_{m+1}^{:,i+1:p})}^{T}}{{\mathbf{A}}^{m}}{{({{\mathbf{A}}^{m}})}^{+}}-{{(\mathbf{H}_{m+1}^{:,i+1:p})}^{T}}{{\mathbf{A}}^{m}}{{({{\mathbf{A}}^{m}})}^{+}}=\mathbf{0}.
\end{equation}
Obviously, (\ref{HAAminusHAA}) deduced from (\ref{CAA_zeros122345}) always holds.
%is equal to zero,
%and then
Thus (\ref{CAA_zeros122345}) has been verified.

After verifying (\ref{CAA_zeros122345}),
let us go back to (\ref{c_wave_defIAACFFCh}) that has been deduced from (\ref{c_col_define4myChol}).
 %we continue to
 %consider
We focus on the entry
\begin{equation}\label{IAACFFC}
{{\mathbf{I}}_{}}-{{\mathbf{A}}^{m}}{{({{\mathbf{A}}^{m}})}^{+}}-{{\mathbf{C}}_{k-1}}{{\mathbf{G}}_{k-1}}\mathbf{G}_{k-1}^{T}\mathbf{C}_{k-1}^{T}
\end{equation}
in (\ref{c_wave_defIAACFFCh}).
% In (\ref{c_wave_defIAACFFCh}),
Substitute (\ref{comments_myB_def111}) into (\ref{IAACFFC}) to obtain
\begin{equation}\label{IAACB2consider123}
{{\mathbf{I}}_{}}-{{\mathbf{A}}^{m}}{{({{\mathbf{A}}^{m}})}^{+}}-{{\mathbf{C}}_{k-1}}\mathbf{B}_{k-1}^{T},
\end{equation}
into which substitute (\ref{C_matrix_def_onlypartH})
%into (\ref{IAACB2consider123})
 to obtain
 \begin{displaymath}
{{\mathbf{I}}_{}}-{{\mathbf{A}}^{m}}{{({{\mathbf{A}}^{m}})}^{+}}-\left( \mathbf{H}_{m+1}^{:,1:k\text{-}1}-{{\mathbf{A}}^{m}}{{\mathbf{D}}_{k\text{-}1}} \right)\mathbf{B}_{k-1}^{T},
 \end{displaymath}
%${{\mathbf{I}}_{}}-{{\mathbf{A}}^{m}}{{({{\mathbf{A}}^{m}})}^{+}}-\left( %\mathbf{H}_{m+1}^{:,1:k\text{-}1}-{{\mathbf{A}}^{m}}{{\mathbf{D}}_{k\text{-}1}} \right)\mathbf{B}_{k-1}^{T}$,
i.e.,
  \begin{equation}\label{IAD21341}
{{\bf{I}}_N} - \left[ {\begin{array}{*{20}{c}}
{{{\bf{A}}^m}}&{{\bf{H}}_{m + 1}^{:,1:k{\rm{ - }}1}}
\end{array}} \right]\left[ \begin{array}{clcr}
{({{\bf{A}}^m})^ + } - {{\bf{D}}_{k - 1}}{\bf{B}}_{k - 1}^T\\
{\bf{B}}_{k - 1}^T
\end{array} \right].
 \end{equation}
%${{\bf{I}}_N} - \left[ {\begin{array}{*{20}{c}}
%{{{\bf{A}}^m}}&{{\bf{H}}_{m + 1}^{:,1:k{\rm{ - }}1}}
%\end{array}} \right]\left[ \begin{array}{l}
%{({{\bf{A}}^m})^ + } - {{\bf{D}}_{k - 1}}{\bf{B}}_{k - 1}^T\\
%{\bf{B}}_{k - 1}^T
%\end{array} \right]$,
 %into which
Then substitute (\ref{BLK_general_inv_def1223}) and (\ref{A_k_def12433})    into (\ref{IAD21341})        to obtain
% \begin{displaymath}
${{\mathbf{I}}_{}}-\mathbf{A}_{k-1}^{m+1}{{(\mathbf{A}_{k-1}^{m+1})}^{+}}$,
% \end{displaymath}
 which is substituted into (\ref{c_wave_defIAACFFCh}) to obtain
${{\mathbf{\tilde{c}}}_{:k}}=\left( {{\mathbf{I}}_{}}-\mathbf{A}_{k-1}^{m+1}{{(\mathbf{A}_{k-1}^{m+1})}^{+}} \right){{\mathbf{h}}_{:k}}$, i.e., (\ref{c_wave_k_final_compute123}).
Since (\ref{c_wave_k_final_compute123}) deduced from (\ref{c_col_define4myChol}) holds,
%\begin{equation}\label{comments_c_waveMine_final_deduce}
%{{\mathbf{\tilde{c}}}_{:k}}={{\mathbf{h}}_{:k}}-\mathbf{A}_{k-1}^{m+1}{{(\mathbf{A}_{k-1}^{m+1})}^{+}}{{\mathbf{h}}_{:k}}.
%\end{equation}
%%Into which
%Finally substitute (\ref{d_col_define4311}) into  (\ref{comments_c_waveMine_final_deduce}) to obtain
%${{\mathbf{\tilde{c}}}_{:k}}={{\mathbf{h}}_{:k}}-\mathbf{A}_{k-1}^{m+1}{{\mathbf{\tilde{d}}}_{:k}}$,
%which is identical to (\ref{c_col_define1}).
%Thus
 we have verified (\ref{c_col_define4myChol}).

%Firstly let us consider $\mathbf{\tilde{c}}_{:k}^{T}{{\mathbf{\tilde{c}}}_{:k}}$ in (\ref{ita_myChol4c_wave}), into which we substitute (\ref{c_col_define4myChol})
%To
Secondly, we
%Now we only need to
 deduce (\ref{c_lentghcCFFcc}).
%we
%Let us
 Substitute (\ref{c_col_define4myChol}) into $\mathbf{\tilde{c}}_{:k}^{T}{{\mathbf{\tilde{c}}}_{:k}}$
 %in (\ref{c_lentghcCFFcc}),
to obtain
$\mathbf{\tilde{c}}_{:k}^{T}{{\mathbf{\tilde{c}}}_{:k}}={{({{\mathbf{c}}_{:k}}-\mathbf{C}_{k-1}^{{}}
{{\mathbf{G}}_{k-1}}\mathbf{G}_{k-1}^{T}\mathbf{C}_{k-1}^{T}{{\mathbf{c}}_{:k}})}^{T}}({{\mathbf{c}}_{:k}}-\mathbf{C}_{k-1}^{{}}{{\mathbf{G}}_{k-1}}\mathbf{G}_{k-1}^{T}\mathbf{C}_{k-1}^{T}{{\mathbf{c}}_{:k}})$, i.e.,
\begin{multline}\label{C_wave_length12312}
\mathbf{\tilde{c}}_{:k}^{T}{{\mathbf{\tilde{c}}}_{:k}}
=\mathbf{c}_{:k}^{T}{{\mathbf{c}}_{:k}}-2\mathbf{c}_{:k}^{T}\mathbf{C}_{k-1}^{{}}{{\mathbf{G}}_{k-1}}
\mathbf{G}_{k-1}^{T}\mathbf{C}_{k-1}^{T}{{\mathbf{c}}_{:k}}+ \\
\mathbf{c}_{:k}^{T}\mathbf{C}_{k-1}^{{}}{{\mathbf{G}}_{k-1}}
\mathbf{G}_{k-1}^{T}\mathbf{C}_{k-1}^{T}\mathbf{C}_{k-1}^{{}}{{\mathbf{G}}_{k-1}}\mathbf{G}_{k-1}^{T}\mathbf{C}_{k-1}^{T}{{\mathbf{c}}_{:k}}.
\end{multline}
Then substitute (\ref{comments_myB_def111}) into the last entry in (\ref{C_wave_length12312}) to write it as
\begin{equation}\label{cCCCFFCc1256}\mathbf{c}_{:k}^{T}\mathbf{C}_{k-1}^{{}}\mathbf{C}_{k-1}^{+}\mathbf{C}_{k-1}^{{}}{{\mathbf{G}}_{k-1}}\mathbf{G}_{k-1}^{T}\mathbf{C}_{k-1}^{T}{{\mathbf{c}}_{:k}}.
\end{equation}
$\mathbf{C}_{k-1}^{+}$ must satisfy (\ref{MoorePenroseDef121a}), i.e.,
%$\mathbf{A}{{\mathbf{A}}^{+ }}\mathbf{A}=\mathbf{A}$,
$\mathbf{C}_{k-1}\mathbf{C}_{k-1}^{+}\mathbf{C}_{k-1}=\mathbf{C}_{k-1}$,
 which can be substituted into (\ref{cCCCFFCc1256}) to simplify  it into
\begin{equation}\label{cCFFCc434316}
\mathbf{c}_{:k}^{T}\mathbf{C}_{k-1}^{{}}{{\mathbf{G}}_{k-1}}\mathbf{G}_{k-1}^{T}\mathbf{C}_{k-1}^{T}{{\mathbf{c}}_{:k}}.
\end{equation}
Finally we can replace the last entry in (\ref{C_wave_length12312}) by (\ref{cCFFCc434316}), to write (\ref{C_wave_length12312}) as
$\mathbf{\tilde{c}}_{:k}^{T}{{\mathbf{\tilde{c}}}_{:k}}=\mathbf{c}_{:k}^{T}
{{\mathbf{c}}_{:k}}-2\mathbf{c}_{:k}^{T}\mathbf{C}_{k-1}^{{}}{{\mathbf{G}}_{k-1}}\mathbf{G}_{k-1}^{T}
\mathbf{C}_{k-1}^{T}{{\mathbf{c}}_{:k}}+\mathbf{c}_{:k}^{T}\mathbf{C}_{k-1}^{{}}{{\mathbf{G}}_{k-1}}\mathbf{G}_{k-1}^{T}\mathbf{C}_{k-1}^{T}{{\mathbf{c}}_{:k}}$, i.e., (\ref{c_lentghcCFFcc}).
\ifCLASSOPTIONcaptionsoff
  \newpage
\fi

\end{document}